\documentclass[a4paper,10pt]{article}
\usepackage[english]{babel}
\usepackage{amsmath, amssymb, amsthm, amscd}
\usepackage[dvips]{graphicx}
\usepackage[active]{srcltx}
\usepackage{pstricks}
\usepackage{pst-all}
\usepackage{amsrefs}
\usepackage{geometry}
\author{DELABARRE Ludovic}

\newcommand{\CQFD}{\hfill{$\Box$}}
\newcommand{\ds}{\displaystyle}
\newcommand{\bs}{\boldsymbol}

\title{On the domain of meromorphy of a multivariate Euler product of Igusa type.}

\newtheorem{prop}{Proposition}
\newtheorem{lemme}{Lemma}
\newtheorem{theo}{Theorem}

\theoremstyle{definition}
\newtheorem{defi}{Definition}
\newtheorem*{dem}{Proof}

\newtheorem{conj}{Conjecture}
\newtheorem{rqs}{Remark}

\theoremstyle{remark}

\newcommand{\Nr}{\mathbb{N}^r \setminus \left\{\textbf{0}\right\}}
\newcommand{\Nn}{\mathbb{N}^{n+1} \setminus \left\{\textbf{0}\right\}}
\newcommand{\summ}{\displaystyle\sum}

\begin{document}

\title{On the domain of meromorphy of a multivariate Euler product of Igusa type.}

\author{Ludovic DELABARRE\footnote{Universit\'{e} de Saint-Etienne, Facult\'{e} des Sciences, Laboratoire de Math{\'e}matiques LaMuse, 23 rue du Docteur Paul Michelon 42023 Sant-Etienne cedex 2, FRANCE. Email: ludovic.delabarre(a)univ-st-etienne.fr 
 \textbf{Keywords:} multivariables Euler products, meromorphic continuation, natural boundary, cyclotomic polynomial. 
\textbf{Classif. math.}: 11M32 11M41 32D15 11N99}}
\date{September 1st, 2011}

\maketitle

\abstract{This work is an answer to a problem posed by N. Kurokawa and H. Ochiai concerning the natural boundary of meromorphy of a multivariate Euler product of Igusa type. More generally, we introduce and determine the maximal domain of meromorphy of a class of multivariate pseudo-uniform Euler products.
}

\tableofcontents

\section{Introduction.}
The aim of this work is to study the domain of meromorphy of some pseudo-uniform Euler products of many variables of the form:
$\ds (s_1,\dots,s_n)\longmapsto \prod_{p \ \textrm{prime}} h(p^{-s_1},\dots,p^{-s_n},p^{-c}),$ where $h(X_1,\dots,X_{n},X_{n+1})\in \mathbb{Z}[X_1,\dots,X_n,X_{n+1}]$ and $c\in \mathbb{Z}\setminus \{0\}$\footnote{Here ``pseudo-uniform'' means that the coefficients of $h$ do not depend on $p$, but the expression $h(p^{-s_1},\dots,p^{-s_n},p^{-c})$ depends both on $p^{-s_i}$ ($i=1,\dots,n$) and on $p$.}.

When $n=1$, the following conjectures formulated by Z. Rudnick and M. du Sautoy (see for example \cite{sautoy3}, 1.4 ; \cite{sautoy}, 1.11) predicts the maximal domain of meromorphy of products of the form
$
Z\left(s\right)=\prod_{p \ \textrm{prime}}h\left(p^{-s},p\right),
$
where $h\left(X_1,X_2\right) \in \mathbb{Z}[X_1,X_2]$.
 
\begin{conj}\label{afconj-rds}
 $ Z\left(s\right)=\prod_{p \ \textrm{prime}}h\left(p^{-s},p\right)$ can be meromorphically continued to the whole complex plane if and only if there exist cyclotomic polynomials  $g_i\left(U\right)$ $\left(i=1,\dots,m\right)$ (meaning divisors of  $\left(1-U^{m_i}\right)^{n_i}$ for a certain ${n_i}$ and a certain $m_i$) and integers $u_i,v_i$ such that:
$
h\left(X_1,X_2\right) = g_1\left(X_1^{u_1}X_2^{v_1}\right)\cdots g_m\left(X_1^{u_m}X_2^{v_m}\right).
$
\end{conj}

\begin{conj}\label{afconj-rds2}
 Let $h(X_1,X_2)=1+\sum_{i=1}^{r}(a_{i0}+a_{i 1}X_1+\cdots + a_{i n_i}X_1^{n_i})X_2^i\in \mathbb{Z}[X_1,X_2]$ which is not a finite product of cyclotomic polynomials as in Conjecture \ref{afconj-rds} and suppose that all cyclotomic factors of $h(X_1,X_2)$ have been removed.
Let $\beta=\textrm{max}\left\{\frac{n_i}{i}: i\in \{1,\dots,r\}\right\}$.
Then $Z\left(s\right)=\prod_{p \ \textrm{prime}}h\left(p^{-s},p\right)$ admits $\Re(s)=\beta$ as natural boundary of meromorphy; meaning that $Z(s)$ can be meromorphically continued to $\{s\in \mathbb{C}: \Re(s)>\beta\}$ but there does not exist any meromorphic continuation beyond the line $\Re(s)=\beta$.
\end{conj}

Eighty years ago, Theodor Estermann (\cite{est}) determined completely the exact domain of meromorphy of the uniform Euler products of one variable $\prod_{p}h(p^{-s})$ ($h(X)\in \mathbb{Z}[X]$): if $h$ is a finite product of cyclotomic polynomials then the corresponding product has a meromorphic continuation to whole $\mathbb{C}$ whereas if $h$ is not then the line $\Re(s)=0$ is a natural boundary of meromorphy.

Many years later, some results, obtained by Bhowmik, Essouabri and Lichtin in \cite{bel} then completed by the author in \cite{delabarre}, have generalized Estermann's theorem to uniform Euler product of many variables.

In this paper, we are interested in a multivariate analogue of the class of products considered in Conjecture \ref{afconj-rds} and \ref{afconj-rds2} and we consider the maximal domain of meromorphy of a pseudo-uniform Euler product of the following form:

$$
Z(s_1,\dots,s_n)=\prod_{p \ \textrm{prime}}h\left(p^{-s_1},\dots,p^{-s_n},p^{-c}\right) := Z^{n+1}(s_1,\dots,s_n,c),
$$ 
where $n>1$, $c \in \mathbb{Z}\setminus \{0\}$ is a fixed nonzero integer and $h(X_1,\dots,X_{n+1}) \in \mathbb{Z}[X_1,\dots,X_{n+1}]$ is a polynomial with integral coefficients of constant coefficient equal to $1$.

What has mainly motivated this study is the resolution of a problem posed by N. Kurokawa and H. Ochiai (see \cite{kurokawa} page 12).

If $A$ is a ring, the multivariate global Igusa zeta function is defined as follows (for $n>1$):
$$
Z^{\textrm{ring}}(s_1,\dots,s_n; A) := \sum_{m_1,\dots,m_n \geq 1} \left |\textrm{Hom}_{\textrm{ring}}\left(A, \frac{\mathbb{Z}}{m_1\cdots m_n \mathbb{Z}}\right)\right | m_1^{-s_1}\cdots m_n^{-s_n}.
$$

By the Chinese remainder theorem, we know that this zeta function can be expressed as an Euler product:
$$
Z^{\textrm{ring}}(s_1,\dots,s_n; A) = \prod_{p}Z_p^{\textrm{ring}}(s_1,\dots,s_n; A)
$$
where 
$$
Z_p^{\textrm{ring}}(s_1,\dots,s_n; A) = \sum_{k_1,\dots,k_n \geq 0} \left |\textrm{Hom}_{\textrm{ring}}\left(A, \frac{\mathbb{Z}}{p^{k_1+\cdots + k_n}\mathbb{Z}}\right)\right|p^{-k_1 s_1 - \cdots - k_n s_n}.
$$

In particular the problem posed in \cite{kurokawa} page $12$ consists in establishing the maximal domain of meromorphy of 
$$
Z^{\textrm{ring}}(s_1,\dots,s_n; \mathbb{Z}[T,T^{-1}]) = \sum_{m_1,\dots,m_n \geq 1} \frac{\varphi(m_1\cdots m_n)}{m_1^{s_1} \cdots m_n^{s_n}},
$$
where $\varphi$ designates the classical Euler function.

As it is pointed out in \cite{kurokawa} (page $12$), the analytic behavior of this product is complicated; nevertheless we manage to determine here its domain of meromorphy by taking full advantage of methods which have been developped in \cite{delabarre} to describe the maximal domain of meromorphy of multivariable uniform Euler products of the form $\prod_{p}h(p^{-s_1},\dots,p^{-s_n})$.

\textbf{Notation:}

In what follows we will use these notations:

For $r\geq 1$ and $n>1$ we write:

$
\displaystyle h(\textbf{X})  =  \displaystyle h(X_1,\dots,X_{n+1})  =  \displaystyle 1+ \sum_{j=1}^{r}a_j X_1^{\alpha_{1, j}}\cdots X_{n+1}^{\alpha_{n+1, j}};
$

with  $\bs{\alpha}_{\cdot j} = (\bs{\alpha}_{(n) j},\alpha_{n+1, j}) = (\alpha_{1, j},\dots,\alpha_{n, j},\alpha_{n+1, j}) \in \mathbb{N}^{n+1}\setminus \{\textbf{0}\}$ for $j \in \{1,\dots, r\}$ and $a_j \in \mathbb{Z}$.

We put $\bs{\alpha}:=\left(\alpha_{\ell,j}\right)_{(\ell,j)\in \{1,\dots,n+1\}\times\{1,\dots,r\}}\in \mathbb{M}_{n+1,r}(\mathbb{N})$ the matrix encoding the exponents of $h$.

For $\ell\in\{1,\dots,n+1\}$, the $\ell$-th row of this matrix will be written $\bs{\alpha}_{\ell \cdot}=(\alpha_{\ell, 1},\dots,\alpha_{\ell, r})$.

For $j\in \{1,\dots,r\}$, the $j$-th column of $\bs{\alpha}$ will be written $\bs{\alpha}_{\cdot j}={}^t \!(\alpha_{1, j},\dots,\alpha_{n, j},\alpha_{n+1, j})={}^t \!({}^t \!\bs{\alpha}_{(n) j}, \alpha_{n+1,j})$, where ${}^t \!\bs{\alpha}_{(n) j}=(\alpha_{1,j},\dots,\alpha_{j,n})$ denotes the first $n$ components of the vector $\bs{\alpha}_{\cdot j}$.
By setting, for $j\in \{1,\dots,r\}$, $\ds\textbf{X}^{\boldsymbol{\alpha}_{\cdot j}} := X_{1}^{\alpha_{1,j}}X_{2}^{\alpha_{2,j}}\cdots X_{n+1}^{\alpha_{n+1,j}}$, we obtain:
\begin{displaymath}
% use packages: array
\begin{array}{lllll}
\displaystyle h(\textbf{X}) & = & \displaystyle h(X_1,\dots,X_{n+1}) & = & \displaystyle 1+ \sum_{j=1}^{r}a_j \textbf{X}^{\bs{\alpha}_{\cdot j}} \\ 
 &  &  & = & \displaystyle 1+\sum_{j=1}^{r}a_j \left(X_1\cdots X_n\right)^{\bs{\alpha}_{(n) j}}X_{n+1}^{\alpha_{n+1, j}}.
\end{array}
\end{displaymath}

For  $\textbf{m}=(m_1,\dots,m_r) \in \mathbb{N}^{r}$, we put $
\Vert \textbf{m} \Vert = \sum_{j=1}^{r}m_{j}.
$

 $\textrm{For} \  \mathbf{s} \in \mathbb{C}^{n+1}, \mathbf{s}= (\mathbf{s}_{(n)},s_{n+1}) = \left(s_1,...,s_{n+1}\right)$ and for  $\ell \in \left\{1,...,n+1\right\}$ we write:

\begin{displaymath}
% use packages: array
\begin{array}{ccccccc}
\sigma_{\ell} & = &\Re\left(s_{\ell}\right); & & \tau_{\ell} &=& \Im\left(s_{\ell}\right);   \\ 
\bs{\sigma} &=& \Re\left( \mathbf{s}\right) &  = & (\bs{\sigma}_{(n)},\sigma_{n+1}) &  = & \left(\sigma_1,\dots,\sigma_{n+1}\right); \\ 
\bs{\tau} &=& \Im\left( \mathbf{s}\right) &  = & (\bs{\tau}_{(n)},\tau_{n+1}) &  = & \left(\tau_1,\dots, \tau_{n+1}\right).
\end{array}
\end{displaymath}

Finally we recall for $\boldsymbol{\nu}=(\nu_1,\dots,\nu_m)$ and $\boldsymbol{w}={}^t \!(w_1,\dots,w_m)$ the classical matrix product between $\boldsymbol{\nu}$ and $\boldsymbol{w}$:
\begin{center}
$\ds \boldsymbol{\nu}\cdot\boldsymbol{w}  := \summ_{i=1}^{m}\nu_{i} w_{i}.$
\end{center}

We must underline the natural appearing of a supplementary hypothesis which permits to distinguish Conjectures \ref{afconj-rds} and \ref{afconj-rds2} from their multivariate analogues since a priori these multivariate analogues contain the conjectures themselves.

So from now on we will suppose that

\begin{equation}\label{afhypnaturelle}
 \textrm{Rank}\left(\bs{\alpha}_{(n) j}, j\in \{1,\dots,r\}\right)>1.
\end{equation}

Indeed, if this hypothesis is not satisfied, we would have the existence of $e$ such that for all $j\in\{1,\dots,r\}$, $\bs{\alpha}_{(n) j} = q_j \bs{\alpha}_{(n) e}$ $(q_j \in \mathbb{Q})$; which would give $
h(\textbf{X})=1+\sum_{j=1}^{r}a_j \left((X_1\cdots X_n)^{\bs{\alpha}_{(n) e}}\right)^{q_j}X_{n+1}^{\alpha_{n+1, j}}.
$
And we would be led to study a one variable product of the form $\prod_{p}h(p^{-s},p^{-c})$.

The aim of this work is to establish the maximal domain of meromorphy of products $\prod_{p}h(p^{-s_1},\dots,p^{-s_n},p^{-c})$ ($n>1$) which cannot be reduced to a one variable product.

$$
\textrm{We also suppose that for all } j\in \{1,\dots,r\}, \bs{\alpha}_{(n) j} = (\alpha_{1, j},\dots,\alpha_{n, j})\neq \textbf{0}.
$$

\begin{rqs}\label{presentationarg}
We must notice that, contrary to the one variable case, the multivariate case permits to take full advantage of methods developped in \cite{delabarre}.
 These methods consist in considering the product in a suitable direction in a neighbourhood of a point of the supposed natural boundary.
The multivariate framework allows to move a point lying on the boundary if necessary and thus it is possible to restrict our attention to ``good'' points (i.e. points for which we are able to prove that no meromorphic continuation is possible beyond) provided that (and it is crucial) these points are generic points (see Definition \ref{defigeneric} and Remark \ref{defigeneric2}) on the boundary.
\end{rqs}

\begin{defi}\label{afcyclo}
 We will say that $h\left(X_1,\dots,X_{n+1}\right)$ is cyclotomic if there exists a finite subset $I$ of $\Nn$ such that we have $
h\left(X_1,\dots,X_{n+1}\right) = \prod_{\lambda=\left(\lambda_1,\dots,\lambda_n\right) \in I}\left(1-X_{1}^{\lambda_1}\cdots X_{n}^{\lambda_{n+1}}\right)^{\gamma\left(\lambda\right)},
$
where the $\gamma\left(\lambda\right)$ for $\lambda \in I$ are positive or negative integers.
\end{defi}

If $h$ is cyclotomic, it is easy to see that the corresponding Euler product is a finite product of classical Riemann zeta function; and consequently it can be meromorphically continued to whole $\mathbb{C}^n$.

So from now on, we will suppose that $h$ \textit{\textbf{is not cyclotomic}} and \textit{\textbf{does not contain any cyclotomic factor}}.

\begin{defi}

For all $\delta \geq 0$ write $
\mathbf{W}(\delta) = \left\{\mathbf{s} \in \mathbb{C}^{n+1}:   \bs{\sigma}\cdot  \bs{\alpha}_{\cdot j} >\delta, \forall  j \in \left\{1,...,r\right\}  \right\};
$
and
$
\mathbf{W}_{c}(\delta) = \left\{\mathbf{s}_{(n)} \in \mathbb{C}^n\mid \bs{\sigma}_{(n)}\cdot  \bs{\alpha}_{(n) j} + c \ \alpha_{n+1, j}>\delta, \forall j \in \{1,\dots,r\} \right\}.
$

\end{defi}

\section{Statements of main results.}

It is straightforward to check that $Z(s_1,\dots,s_n)$ is holomorphic on $\mathbf{W}_c(1)$.
Moreover, we will prove (see Theorem \ref{afana}) that one can continue  meromorphically  $Z(s_1,\dots,s_n)$ to $\mathbf{W}_c(0)$.
The main result of this paper consists in verifying, in most cases, that the edge $\partial \mathbf{W}_c(0)$ of $\mathbf{W}_c(0)$ is a natural boundary of meromorphy.

Before announcing the results, firstly we introduce a definition.

Since $\mathbf{W}(0) = \{\mathbf{s}\in \mathbb{C}^{n+1} : \Re(\mathbf{s}\cdot  \bs{\alpha}_{\cdot j})\geq 0, \forall j=1,\dots,r\}$, then the edge $\partial \mathbf{W}(0)$ of $\mathbf{W}(0)$ is a polyhedron whose faces are of the form $\ds
\mathcal{F}(\bs{\alpha}_{\cdot e}) = \{\mathbf{s}\in \overline{\mathbf{W}(0)}:\Re(\mathbf{s}\cdot  \bs{\alpha}_{\cdot e})=0\};
$
for a vector $\bs{\alpha}_{\cdot e}\in \{\bs{\alpha}_{\cdot 1},\dots,\bs{\alpha}_{\cdot r}\}$.
We will say by abuse of language that $\mathcal{F}(\bs{\alpha}_{\cdot e})$ is a face of polar vector $\bs{\alpha}_{\cdot e}$.

 Now let $\mathcal{F}(\bs{\alpha}_{\cdot e})$ be a face of the edge $\partial \mathbf{W}(0)$ of $\mathbf{W}(0)$ as above and consider in particular $\widehat{\bs{\alpha}}_{\cdot e}\in \mathbb{N}^n, \widehat{\bs{\alpha}}_{\cdot e} \in \mathbb{Q}\bs{\alpha}_{\cdot e}$ the vector collinear with  $\bs{\alpha}_{\cdot e}$ whose nonzero components are relatively prime.

\begin{defi}\label{afe-iemepartie}
 Given $e \in \{1,\dots,r\}$ we denote by $\langle\bs{\alpha}_{\cdot e}\rangle$ the line connecting $\mathbf{0}$ and the
integer point $\bs{\alpha}_{\cdot e}$ in $\mathbb{R}^n$, and then define the $e$-th main part of $h$ 
$\ds
[h]_{e}(\textbf{X}) = \sum_{\bs{\alpha}_{\cdot j} \in \langle\bs{\alpha}_{\cdot e}\rangle} a_j \textbf{X}^{\bs{\alpha}_{\cdot j}}.
$
\end{defi}

\begin{defi}
Given $e \in \{1,\dots, r\}$ we set  
\begin{displaymath}
% use packages: array
\begin{array}{ll}
\displaystyle  \Lambda_{e}  &= \big\{j \in \left\{1,\dots,r\right\} : \bs{\alpha}_{\cdot j} \in  \langle\bs{\alpha}_{\cdot e}\rangle \big\}\\
B_e &= \big\{\bs{\beta} \in \mathbb{N}^r : \beta_j = 0 \quad  \text{if} \quad j \notin \Lambda_e \big\}.
\end{array}
\end{displaymath}
\end{defi}

It is clear that for all $j\in \Lambda_e$ there exists $q_j\in \mathbb{N}^{*}$ such that $\bs{\alpha}_{\cdot j}=q_j \widehat{\bs{\alpha}}_{\cdot e}$.
Then we define
$
\ds \widetilde{[h]_e}(T)  :=  \ds 1+\sum_{j\in \Lambda_e}a_j T^{q_j}\in \mathbb{Z}[T] \ \textrm{verifying} \ \widetilde{[h]_e}(\textbf{X}^{\widehat{\bs{\alpha}}_{\cdot e}}) = [h]_e(\textbf{X}).
$

\begin{defi}\label{anouvcond2}

We will say that the face $\mathcal{F}(\bs{\alpha}_{\cdot e})$ is a non-degenerate face if the one variable polynomial $\widetilde{[h]_e}(T)$ has no multiple root.
 \end{defi}

We suppose here that $h$ is such that $\partial \mathbf{W}(0)$ contains at least one non-degenerate face $\mathcal{F}(\bs{\alpha}_{\cdot e})$  in the sense of Definition \ref{anouvcond2}.

The aim of this article is essentially to prove two complementary results concerning the natural boundary (meaning the boundary beyond which there does not exist any meromorphic continuation) of $\prod_{p}h(p^{-s_1},\dots,p^{-s_n},p^{-c})$ which depend on the validation of an hypothesis that we will note (H) (see Theorem \ref{afprincipal}).

 We will see that if this property (H) is satisfied we are able to determine the natural boundary in a strong sense (see Theorem \ref{afprincipal}) whereas if it is not verified, we still obtain the natural boundary but in a weaker sense (see Theorem \ref{afprincipal2}): we will see that it cannot exist any meromorphic continuation by translating the boundary to the left.

\begin{theo}\label{afprincipal}
 Let $ c \in \mathbb{Z}\setminus \{0\}$ and  $
Z(s_1,\dots,s_n) = \prod_{p}h(p^{-s_1},\dots,p^{-s_n},p^{-c}).
$
The product $(s_1,\dots,s_n)\longmapsto Z(s_1,\dots,s_n)$ converges absolutely in the domain $\mathbf{W}_c(1)$ and admits a meromorphic continuation to $\mathbf{W}_c(0)$.
 Moreover, assume that the polynomial $ h(X_1,\dots,X_n,X_{n+1})$ is not cyclotomic, does not contain any cyclotomic factors, admits at least one non-degenerate face $\mathcal{F}(\bs{\alpha}_{\cdot e})$, verifies (\ref{afhypnaturelle}) and satisfies in addition the following property (H):
$$
 \textrm{for all} \ j\in \{1,\dots,r\} \ \textrm{such that} \ \bs{\alpha}_{\cdot j} \notin \mathbb{Q}\bs{\alpha}_{\cdot e}, \  \bs{\alpha}_{(n) j} \notin \mathbb{Q}\bs{\alpha}_{(n) e}.
$$
Then the set $\{(s_1,\dots,s_n)\in \mathbb{C}^n: (s_1,\dots,s_n,c)\in \mathcal{F}(\bs{\alpha}_{\cdot e})\}\subseteq\partial \mathbf{W}_c(0) $ is a natural boundary (in the strong sense): there does not exist any continuation of $Z(s_1,\dots,s_n)$ to a domain containing an open ball $\mathcal{B}$ (of dimension $n$) centered in a point $\mathbf{s}_{(n)}^0$ such that $(\mathbf{s}_{(n)}^0,c)\in \mathcal{F}(\bs{\alpha}_{\cdot e})$.
\end{theo}
%\vspace{0.5cm}

\begin{theo}\label{afprincipal2}
 Let $ c \in \mathbb{Z}\setminus \{0\}$ and $
Z(s_1,\dots,s_n) = \prod_{p}h(p^{-s_1},\dots,p^{-s_n},p^{-c}).
$
The product $(s_1,\dots,s_n)\longmapsto Z(s_1,\dots,s_n)$ converges absolutely in the domain $\mathbf{W}_c(1)$ and admits a meromorphic continuation to $\mathbf{W}_c(0)$.
 Moreover, assume that the polynomial $ h(X_1,\dots,X_n,X_{n+1})$ is not cyclotomic, does not contain any cyclotomic factors, verifies (\ref{afhypnaturelle}) and admits at least one non-degenerate face $\mathcal{F}(\bs{\alpha}_{\cdot e})$ but does not satisfy the property (H) of Theorem \ref{afprincipal}.
Suppose in addition the following property:

\begin{footnotesize}\begin{equation}\label{aeglobal}
\textrm{if} \ \bs{\alpha}_{(n) j_0}\notin \mathbb{Q}\bs{\alpha}_{(n) e} \  \textrm{then the polynomials} \ 1+\sum_{\bs{\alpha}_{\cdot j}\in \mathbb{Q}\bs{\alpha}_{\cdot e}}a_j  \textbf{X}^{\bs{\alpha}_{\cdot j}} \ \textrm{and} \ \sum_{j:\bs{\alpha}_{\cdot j}-\bs{\alpha}_{\cdot j_0}\in \mathbb{Q}\bs{\alpha}_{\cdot e}}a_j \textbf{X}^{\bs{\alpha}_{\cdot j}} \  \textrm{are relatively prime.}
\end{equation}              \end{footnotesize}

Then the edge $\partial \mathbf{W}_c(0)$ of $\mathbf{W}_c(0)$ is a natural boundary (in the weak sense): $Z(s_1,\dots,s_n)$ does not admit a meromorphic extension to $\mathbf{W}_{c}(\delta)$ for any $\delta<0$.
In particular, $Z(s_1,\dots,s_n)$ does not admit any meromorphic continuation to $\mathbb{C}^n$.
\end{theo}

\vspace{0,4cm}

As an application, we will see that we can determine the natural boundary (in the strong sense) of Igusa's zeta function $Z^{\textrm{ring}}(s_1,\dots,s_n; \mathbb{Z}[T,T^{-1}])$ by  obtaining the following result:

\begin{theo}\label{afapplication}
 The maximal domain of meromorphy $\mathcal{M}$ of Igusa's zeta function:
$$
Z^{\textrm{ring}}(s_1,\dots,s_n; \mathbb{Z}[T,T^{-1}]) = \sum_{m_1,\dots,m_n \geq 1} \frac{\varphi(m_1\cdots m_n)}{m_1^{s_1} \cdots m_n^{s_n}}
$$
is given by $\ds
\mathcal{M} = \left\{(s_1,\dots,s_n)\in \mathbb{C}^n\mid \forall \ I\subseteq \{1,\dots,n\}, \sum_{i \in I}\sigma_i > -1+\#I\right\}.
$
In particular, if $(s_1^0,\dots,s_n^0) \in \partial \mathcal{M},$ then there does not exist any meromorphic continuation of $Z^{\textrm{ring}}(s_1,\dots,s_n; \mathbb{Z}[T,T^{-1}])$ to a domain containing an open ball $\mathcal{B}$ of dimension $n$ centered in $(s_1^0,\dots,s_n^0)$.
\end{theo}

\section{Proof of Theorem \ref{afprincipal}.}\label{afproofthm1}
\subsection{Meromorphic continuation of $Z(s_1,\dots,s_n)$.}

The pseudo-uniform Euler product $(s_1,\dots,s_n)\longmapsto Z(s_1,\dots,s_n)$, absolutely convergent in $\mathbf{W}_c(1)$, admits a meromorphic continuation to $\mathbf{W}_c(0)$. This fact follows directly from the meromorphic continuation of the uniform product $(s_1,\dots,s_n,s_{n+1}) \longmapsto \prod_{p}h(p^{-s_1},\dots,p^{-s_n},p^{-s_{n+1}})$ to $\mathbf{W}(0)$ as it is proved in \cite{delabarre} (see \cite{delabarre}, Theorem 1).
Furthermore, we are able to give an expression of the meromorphic extension of  $Z^{n+1}(s_1,\dots,s_{n+1})$ in  $\mathbf{W}(\delta)$ for all $\delta > 0$, and hence an expression of the continuation of $Z(s_1,\dots,s_n)$ in $W_c(\delta)$ for all $\delta>0$.

The key point of the continuation of $Z^{n+1}$ is the writing of $h(\textbf{X})$ as an infinite product a cyclotomic factors:

\begin{prop}\label{expansion}
Consider the quantity:
\begin{equation}\label{afdist-h-zero}
C:=C(h) = \frac{1}{|a_1|+\cdots + |a_r|}.
\end{equation}
 If  $|\textbf{X}^{\bs{\alpha}_{\cdot j}}| < C$ for $j \in \{1,\dots,r\}$, then:
\begin{displaymath}
\begin{array}{lll}
h(\textbf{X})=1 + a_1 \textbf{X}^{\bs{\alpha}_{\cdot 1}}+ \dots + a_r \textbf{X}^{\bs{\alpha}_{\cdot r}} &=&\displaystyle \prod_{\bs{\beta}=(\beta_1,\dots,\beta_r) \in \Nr}\left(1-\prod_{\ell=1}^{n+1}X_{\ell}^{\bs{\alpha}_{\ell \cdot}\cdot {}^t \!\bs{\beta}}\right)^{\gamma(\bs{\beta})} \\
& = &\displaystyle \prod_{\bs{\beta}\in \Nr}\left(1-\textbf{X}^{\bs{\alpha}\cdot {}^t \!\bs{\beta}}\right)^{\gamma(\bs{\beta})};
\end{array}
\end{displaymath}
where the right side converges absolutely and
$$
\gamma(\bs{\beta}) = \displaystyle \sum_{\substack{  
\textbf{b} \in \mathbb{N}^r\setminus \left\{\textbf{0}\right\} \\ m \in \mathbb{N} \\ m \textbf{b} = \bs{\beta} }} \left(\left(-1\right)^{\Vert \textbf{b} \Vert } \frac{\mu\left(m\right)}{m} \frac{\left(\Vert \textbf{b} \Vert -1\right)!}{b_{1}!...b_{r}!} a_{1}^{b_{1}}...a_{r}^{b_{r}} \right)\in \mathbb{Z}.
$$
\end{prop}

\begin{dem}
 See \cite{delabarre}, Lemma 2 and Corollary 2.2.
\end{dem}

\begin{theo}\label{afana}

$Z^{n+1}\left(s_1,\dots,s_{n+1}\right)$ is meromorphic in $\mathbf{W}(0)$.

Moreover if we write for all $\delta>0$ $M_\delta=\left[C^{-\frac{1}{\delta}}\right]+1$ $(M_\delta \in \mathbb{N})$, the following relation holds in $\mathbf{W}(\delta)$:
$$
Z^{n+1}\left(s_1,\dots,s_{n+1}\right) = \prod_{p \leq M_\delta}h\left(p^{-s_1},...,p^{-s_{n+1}}\right)\prod_{\bs{\beta}\in \Nr} \zeta_{M_\delta} \left(\mathbf{s}\cdot \bs{\alpha}\cdot {}^t \! \bs{\beta} \right)^{-\gamma(\bs{\beta})};
$$

where $\zeta_{M_{\delta}}(z) = \zeta(z)\prod_{p\leq M_{\delta}}\left(1-p^{-z}\right)$ ($\zeta$ being the classical Riemann zeta function) has exactly the same zeros and poles as the classical Riemann zeta function with the same multiplicities.
In addition, the possible zeros or poles of $\prod_{\bs{\beta}\in \Nr}\zeta_{M_\delta} \left(\mathbf{s}\cdot \bs{\alpha}\cdot {}^t \! \bs{\beta} \right)^{-\gamma(\bs{\beta})}$, which is meromorphic in $\mathbf{W}(\delta)$, belong to the set:
$$
\Phi_{\delta} = \left\{\mathbf{s} \in \mathbf{W}(\delta), \exists \bs{\beta} \in \mathbb{N}^r, \mathbf{s}\cdot \bs{\alpha}\cdot {}^t \! \bs{\beta} = \rho, \ \rho \ \textrm{zero \ or \ pole \ of \ } \zeta\left(.\right) \right\}.
$$ 

\end{theo}

\begin{dem}
 See \cite{delabarre}, Theorem 1.
\end{dem}

\subsection{Determination of the natural boundary of $Z(s_1,\dots,s_n)$.}
Consider a point  $\mathbf{s}^0=(\mathbf{s}_{(n)}^0,c)$ with $\mathbf{s}_{(n)}^0\in \partial \mathbf{W}_c(0)$ lying on a non-degenerate face $\mathcal{F}(\bs{\alpha}_{\cdot e})$ of $\partial \mathbf{W}(0)$ (which exists by hypothesis)  of real part $\bs{\sigma}^0 = (\bs{\sigma}_{(n)}^0,c)$ and of imaginary part $\bs{\tau}^0 = (\bs{\tau}_{(n)}^0,0)$.
Consequently, we have for all $ j \in \left\{1,...,r\right\}$,  $\bs{\sigma}^0\cdot \bs{\alpha}_{\cdot j} \geq 0 $ and $\bs{\sigma}^0\cdot \bs{\alpha}_{\cdot e}=0$.
Consider an open ball $\mathcal{B}$ of dimension $n$ and of arbitrarily small radius around the point $\mathbf{s}_{(n)}^0$.
%EVENTUELLEMENT HYP...
Let us start by giving the definition of the notion of generic set which is fundamental in whole this paper:
\begin{defi}\label{defigeneric}
 A subset $G$ of a set $E$ is said to be generic in $E$ if the complement of $G$ in $E$ has empty interior.
\end{defi}

\begin{rqs}\label{defigeneric2}
 In the following when we will use the term ``generic point'' in a set $E$, it will be understood that we consider any point belonging to some generic set in $E$.
\end{rqs}

To begin with, let us check that it is possible we can choose  $\mathbf{s}_{(n)}^0$ in a generic set inside $\mathcal{B}$ such that $\mathbf{s}^0 = (s^0_1,\dots,s^0_n,c) = (\mathbf{s}_{(n)}^0,c) \in \partial \mathbf{W}(0)$ in a way to have:
\begin{equation}\label{afhyp00}
  (\bs{\sigma}^0_{(n)},c)\cdot \bs{\alpha}_{\cdot j} = 0  \ \Longleftrightarrow \  \bs{\alpha}_{\cdot j} \in \mathbb{Q} \bs{\alpha}_{\cdot e}.
\end{equation}

For this assume that we have the existence of $j_1$ and of $j_2$ such that $\mathbb{Q}\bs{\alpha}_{\cdot j_1} \neq \mathbb{Q}\bs{\alpha}_{\cdot j_2}$ and such that
$
 (\bs{\sigma}^0_{(n)},c)\cdot  \bs{\alpha}_{\cdot j_1} =  (\bs{\sigma}^0_{(n)},c)\cdot  \bs{\alpha}_{\cdot j_2} = 0.
$
Then we have:
\begin{enumerate}
 \item The equality $ (\bs{\sigma}^0_{(n)},c)\cdot  \bs{\alpha}_{\cdot j_1} = \summ_{\ell=1}^{n}\sigma^0_{\ell} \alpha_{\ell, j_1} + c \ \alpha_{n+1, j_1} = 0$  defines an affine real space $A_1$ of dimension $n-1>0$ according to $\bs{\sigma}^0_{(n)}$.
 \item The equality $ (\bs{\sigma}^0_{(n)},c)\cdot  \bs{\alpha}_{\cdot j_2} = \summ_{\ell=1}^{n}\sigma^0_{\ell} \alpha_{\ell, j_2} + c \ \alpha_{n+1, j_2} = 0$ defines an affine real space $A_2 \neq A_1$ of dimension $n-1>0$ according to $\bs{\sigma}^0_{(n)}$ because if $\bs{\alpha}_{(n) j_2}\in \mathbb{Q}\bs{\alpha}_{(n) j_1}$ we would have necessarily $\bs{\alpha}_{\cdot j_2}\in \mathbb{Q}\bs{\alpha}_{\cdot j_1}$ since $c\neq 0$.
\end{enumerate}
Consequently we have necessarily $\bs{\sigma}^0_{(n)} \in A_1\bigcap A_2$; and hence it belongs to an affine subspace of dimension less than or equal to $n-2$ of empty interior in $\mathbb{R}^{n-1}$; so we have $(\ref{afhyp00})$ by choosing $\bs{\sigma}^0_{(n)}$ in a generic set such that $\mathbf{s}_{(n)}^0 \in \partial \mathbf{W}_c(0)\cap \mathcal{B}$.

In the same way, we can suppose  (by moving $\bs{\sigma}^0_{(n)}$ if necessary such that $(\bs{\sigma}^0_{(n)},c)\in \partial \mathbf{W}(0)$ by avoiding a countable union of closed sets of empty interior which remains of empty interior according to Baire's theorem) that for $\bs{\lambda}=(\lambda_1,\dots,\lambda_{n+1})\in\mathbb{Q}^{n+1}$:
\begin{equation}\label{afdcond}
 (\bs{\sigma}^0_{(n)},c)\cdot {}^t \! \bs{\lambda}=0 \ \Longleftrightarrow \ \bs{\lambda} \in \mathbb{Q}\bs{\alpha}_{\cdot e};
\end{equation}

\vspace{0.2cm}

\underline{\textbf{Note:}} the previous argument is the simplest illustration of a recurring principle appearing throughout this paper;  a principle briefly presented above (see Remark \ref{presentationarg} page \pageref{presentationarg}) which consists in  moving the point $\mathbf{s}_{(n)}^0$ on the boundary by avoiding a ``bad'' set. 
And the crucial point is to verify that this ``bad'' set has an empty interior in order to prove that the ``good'' points are generic in the set $\{(s_1,\dots,s_n)\in \mathbb{C}^n: (s_1,\dots,s_n,c)\in \mathcal{F}(\bs{\alpha}_{\cdot e})\}\subseteq\partial \mathbf{W}_c(0)$ so that \textit{all} this set $\{(s_1,\dots,s_n)\in \mathbb{C}^n: (s_1,\dots,s_n,c)\in \mathcal{F}(\bs{\alpha}_{\cdot e})\}$ is a natural boundary of meromorphy.

\vspace{0.2cm}

From now on, we will also suppose without loss of generality (rearranging the indexes if necessary) that $\alpha_{n, e} \neq 0$. 

We fix a direction $\bs{\theta}=(\theta_1,\dots,\theta_n,0):=(\bs{\theta}_{(n)},0) \in \mathbb{Q}^{n+1}$ with $\theta_{n+1}=0$  and consider the one variable function in $t$ complex:
$$
t \longmapsto Z^{n+1}(s^0_{1}+t\theta_1,\dots, s^0_{n}+t\theta_n,c) = Z^{n+1}(\mathbf{s}^0+t\bs{\theta}) = Z(s^0_{1}+t\theta_1,\dots, s^0_{n}+t\theta_n) =  Z(\mathbf{s}_{(n)}^0+t\bs{\theta}_{(n)});
$$
for $t$ lying inside a rectangle (for $u \in \mathbb{R}, \eta>0$ ):
\begin{center}
% use packages: array
\begin{tabular}{ll}
$\Xi_{u,\eta}:$ & $0<\Re\left(t\right)< 1 $\\ 
 & $0<u<\Im\left(t\right)<u+\eta.$
\end{tabular}
\end{center}

We suppose that $\bs{\theta}$ satisfies the following conditions:

 \begin{equation}\label{afhyp-theta1}
  \bs{\theta}\cdot  \bs{\alpha}_{\cdot j} = \bs{\theta}_{(n)}\cdot  \bs{\alpha}_{(n) j} \geq 1 \ \textrm{for all} \ j\in\{1,\dots,r\}.
 \end{equation}

We assume also that

\begin{equation}\label{entier}
 \bs{\theta}_{(n)}\cdot  \widehat{\bs{\alpha}}_{(n) e}\in \mathbb{N}^*;
\end{equation}
(where $\widehat{\bs{\alpha}}_{\cdot e}$ is the vector collinear with $\bs{\alpha}_{\cdot e}$ whose components are relatively prime).

In this way, since for $j\in \Lambda_e$ we have seen (see Definition \ref{afe-iemepartie}) that there exists $q_j\in\mathbb{N}^*$ such that $\bs{\alpha}_{\cdot j}=q_j\widehat{\bs{\alpha}}_{\cdot e}$, we have necessarily for all $j\in \Lambda_e$ $\bs{\theta}_{(n)}\cdot  \bs{\alpha}_{(n) j}\in \mathbb{N}^*$.

 \vspace{0.5cm}

The aim of what follows is to prove the existence of an accumulation of zeros or poles of the one variable function  $t\longmapsto Z(\mathbf{s}_{(n)}^0+t\bs{\theta}_{(n)})$  inside the bounded rectangle $\Xi_{u,\eta}$ for any $u\in \mathbb{R}, \eta>0$. 

By putting $s_{n+1}=c$ in Theorem \ref{afana},  we have an expression for $Z(s_1,\dots,s_n)$ inside $\mathbf{W}_c(\delta)$ for all $\delta>0$.
According to the condition (\ref{afhyp-theta1}) satisfied by $\bs{\theta}$, this expression remains well-defined for all $\delta>0$ for $\Re(t)>\delta$ if we put $\mathbf{s}_{(n)}:=\mathbf{s}^0_{(n)}+t\bs{\theta}_{(n)}$  since $\Re(t)>\delta$ implies $\mathbf{s}^0+t\bs{\theta}=(\mathbf{s}^0_{(n)},c)+t(\bs{\theta}_{(n)},0)\in \mathbf{W}_c(\delta)$. This allows to characterize the eventual zeros or poles of $t\longmapsto Z(\mathbf{s}_{(n)}^0+t\bs{\theta}_{(n)})$. Indeed the poles necessarily come from $\zeta$-factors and hence belong to $\Phi_{\delta}$ (with $s_{n+1}=c$) for some $\delta>0$ whereas the zeros of $t\longmapsto Z(\mathbf{s}_{(n)}^0+t\bs{\theta}_{(n)})$ are provided by the $\zeta$-factors and by the zeros  of  $t\longmapsto h(p^{-s^0_1-t\theta_1},\dots,p^{-s^0_n-t\theta_n},p^{-c})$.

For the sake of presentation, consider the following triplet of parameters $$\bs{\mu}=(p,\bs{\tau}_{(n)}^0, \bs{\theta}_{(n)}).$$ When this parameter will be put in index, this will mean a dependence according to the prime number $p$, the imaginary part $\bs{\tau}_{(n)}^0$ of $\mathbf{s}^0_{(n)}$ and the direction $\bs{\theta}_{(n)}$.

\begin{defi}
We define the generalized polynomial $W_{\bs{\mu},\bs{\sigma}^0}(X,Y)$ of $h$ depending on $\mathbf{s}_{(n)}^0, \bs{\theta}_{(n)}$ and $p$ for $X\in \mathbb{C}\setminus \mathbb{R}_{-}$ and $Y\in \mathbb{C}$ deprived of an half-line (i.e. $Y\in \mathbb{C}\setminus e^{i\mathfrak{b}}\mathbb{R}_{+}$ for some $\mathfrak{b}\in [0,2\pi[$)  as
$
\ds W_{\bs{\mu},\bs{\sigma}^0}(X,Y) = 1+\sum_{j=1}^{r}a_j p^{-i\bs{\tau}_{(n)}^0\cdot  \bs{\alpha}_{(n) j}}X^{\bs{\sigma}^0_{(n)}\cdot  \bs{\alpha}_{(n) j}+c\alpha_{n+1, j}}Y^{\bs{\theta}_{(n)}\cdot  \bs{\alpha}_{(n) j}}.
$
\end{defi}
Since $\mathbf{s}^0 = \left(\bs{\sigma}^0_{(n)},c\right) + i \left(\bs{\tau}_{(n)}^0,0\right)$, then for all $p$ prime and $t\in \Xi_{u,\eta}$ such that $\Im(t)\log(p)+\mathfrak{b}\notin 2\pi\mathbb{Z}$ we have
$
\ds W_{\bs{\mu},\bs{\sigma}^0}(p^{-1},p^{-t}) = h\left(p^{-s^0_1-t\theta_1},\dots,p^{-s^0_n-t\theta_n},p^{-c}\right).
$

\vspace{0.3cm}

The aim of what follows is to characterize the zeros of these generalized polynomials $W_{\bs{\mu},\bs{\sigma}^0}(X,Y)=0$ by expressing $Y=\Omega(X)$ as a function of $X$ such that $W_{\bs{\mu},\bs{\sigma}^0}(X,\Omega(X))=0$. The problem is that we cannot apply the classical Puiseux theorem to find these solutions because $W_{\bs{\mu},\bs{\sigma}^0}(X,Y)$ is not a real polynomial and in particular is not well-defined for $X$ in a neighbourhood of $0$.
This is the reason why here we need the hypothesis that the face $\mathcal{F}(\bs{\alpha}_{\cdot e})$ is non-degenerate; which allows us to prove the existence of the solutions and equally to have a good control on the convergence (particularly on the dependence on $p$). This result is given by the following proposition:

\begin{prop}\label{athpuiseuxpolgen} {\bf (Puiseux theorem for $W_{\bs{\mu},\bs{\sigma}^0}\left(X,Y\right)$)}.\\
Let $q\in \mathbb{N}^*$ be the smallest positive integer verifying $q \bs{\theta}_{(n)} \cdot \bs{\alpha}_{(n) j} \in \mathbb{N}^*$ for all $j=1,\dots,r$.
 Consider the finite set: 
$$
\mathfrak{p}_{\bs{\mu}}:=\left\{c_{\bs{\mu}} \in \mathbb{C};~ \exists \mathfrak{r} \ \textrm{root of} \ \widetilde{[h]_e}(T) \ \textrm{such that} 
\ c_{\bs{\mu}}^{q \bs{\theta}_{(n)}\cdot  \widehat{\bs{\alpha}}_{(n) e}} = \mathfrak{r} p^{i\bs{\tau}_{(n)}^0\cdot  \widehat{\bs{\alpha}}_{(n) e}}   \right\};
$$
(where $\widehat{\bs{\alpha}}_{\cdot e}$ is the vector collinear with $\bs{\alpha}_{\cdot e}$ whose components are relatively prime).

There exists $\epsilon_1>0$ (not depending on $p$ nor on $\bs{\tau}^0_{(n)}$) such that  for all $X\in \mathcal{H}:=\{X\in \mathbb{C}\setminus \mathbb{R}_{-}, |X|<\epsilon_1\}$ the equation $W_{\bs{\mu},\bs{\sigma}^0}(X,Y)=0$ admits the set of solutions 
$\ds
Y = \Omega_{\bs{\mu},c_{\bs{\mu}},\bs{\sigma}^0}( X) \ (c_{\bs{\mu}} \in \mathfrak{p}_{\bs{\mu}});
$
where for all $ c_{\bs{\mu}} \in \mathfrak{p}_{\bs{\mu}}$, $X\longmapsto \Omega_{\bs{\mu},c_{\bs{\mu}},\bs{\sigma}^0}( X)$ is an holomorphic function on  $\mathcal{H}$ and satisfies for all  $X\in \mathcal{H}:$ $
\Omega_{\bs{\mu},c_{\bs{\mu}},\bs{\sigma}^0}(X) = \sum_{k = 0}^{\kappa_{c_{\bs{\mu}}}}\mathfrak{c}_{k}(c_{\bs{\mu}},\bs{\mu})X^{\vartheta(\bs{\sigma}^0_{(n)})_k}
$
with:
\begin{enumerate}
 \item $\kappa_{c_{\bs{\mu}}}\in \mathbb{N}\cup\{+\infty\}$;
\item $\vartheta(\bs{\sigma}^0_{(n)})_0=0<\vartheta(\bs{\sigma}^0_{(n)})_1<\cdots$ is a stricly increasing sequence (independent of $\bs{\tau}_{(n)}^0$);
\item $\lim_{k\to +\infty} \vartheta_k = +\infty$ if $\kappa_{c_{\bs{\mu}}}=+\infty$;
\item there exist two constants $D_{\epsilon_0} >1$ and $A(\bs{\sigma}^0_{(n)}) >0$ (independent of $p$, $\bs{\tau}^0_{(n)}$ and $k$) such that 
$
\left|\mathfrak{c}_k (c_{\bs{\mu}},\bs{\mu})\right|\ll  D_{\epsilon_0}^{A(\bs{\sigma}^0_{(n)})\vartheta(\bs{\sigma}^0_{(n)})_k}
$
uniformly in $p$ prime and in $k$;
\item $\mathfrak{c}_{0}(c_{\bs{\mu}},\bs{\mu}) = c_{\bs{\mu}}^q$, in particular
$\ds\left|\mathfrak{c}_{0}(c_{\bs{\mu}},\bs{\mu})\right| = |c|^{\frac{1}{\bs{\theta}_{(n)}\cdot \bs{\widehat{\alpha}}_{(n) e}}}.
$
\end{enumerate}
Moreover $\ds\{ \mathfrak{c}_{0}(c_{\bs{\mu}},\bs{\mu}); c_{\bs{\mu}} \in \mathfrak{p}_{\bs{\mu}} \}
=\{ u \in \mathbb{C}; ~ \exists \mathfrak{r} \ \textrm{root of} \ \widetilde{[h]_e}(T) \ \textrm{such that}
~ u^{q \bs{\theta}_{(n)}\cdot \bs{\widehat{\alpha}}_{(n) e}} = \mathfrak{r}p^{i\bs{\tau}_{(n)}^0\cdot \bs{\widehat{\alpha}}_{(n) e}} \}.$

\end{prop}

\begin{dem}
 See \cite{delabarre}, Proposition 2. \CQFD
\end{dem}

Now the aim is to find a Puiseux series  $Y=\Omega_{\bs{\mu},\bs{\sigma}^0}(X)$ in a neighbourhood of $X=0$ ($X\in \mathbb{C}\setminus \mathbb{R}_{-}$) such that  $W_{\bs{\mu},\bs{\sigma}^0}(X,\Omega_{\bs{\mu},\bs{\sigma}^0}(X))=0$ and verifying
$\ds
|\Omega_{\bs{\mu},\bs{\sigma}^0}(X)|<1 \ \textrm{for} \ X>0 \ \textrm{small enough}.
$

In this way we will have infinitely many zeros $t_{m,\bs{\mu},\bs{\sigma}^0}$ ($m\in\mathbb{Z}$, $p$ prime large enough) of the form
\begin{equation}\label{afzerotmp}
t_{m,\bs{\mu},\bs{\sigma}^0} = -\frac{\log\left(\Omega_{\bs{\mu},\bs{\sigma}^0}(p^{-1})\right)}{\log(p)}+\frac{2 \pi m i}{\log(p)}
\end{equation}

of strictly positive real part inside $\Xi_{u,\eta}$ for $p$ large enough.

So now consider any solution $\Omega_{\bs{\mu},\bs{\sigma}^0}(X)$ of $W_{\bs{\mu},\bs{\sigma}^0}(X,Y)=0$ (in finite number) that we will write as follows:

\begin{equation}\label{afbranchepuiseux}
  \Omega_{\bs{\mu},\bs{\sigma}^0}\left(X\right)=c_{\bs{\mu},0}+c_{\bs{\mu},1} X^{\vartheta_{1}}+\cdots+c_{\bs{\mu},N} X^{\vartheta_{N}}+\Omega_{\bs{\mu},\bs{\sigma}^0,N+1}\left(X\right),    \ \left(N \geq 1\right)
\end{equation}
where $c_{\bs{\mu},m} \in \mathbb{C}$;   $\vartheta_{N}= \vartheta(\bs{\sigma}^0_{(n)})_{N}>\cdots>\vartheta_{1} = \vartheta(\bs{\sigma}^0_{(n)})_{1} \in \mathbb{N}^*$;  $\Omega_{\bs{\mu},\bs{\sigma}^0,N+1}\left(X\right) = o\left(X^{\vartheta_{N}}\right)$;
and we have $\ds
\forall k \in \left\{1,\dots,f\right\}, W_{\bs{\mu},\bs{\sigma}^0}\left(X,\Omega_{\bs{\mu},\bs{\sigma}^0}\left(X\right)\right) =  0.
$

We have to notice that, according to Proposition \ref{athpuiseuxpolgen} (claim 5.), the main term $c_{\bs{\mu},0}$ of a Puiseux branch is a root of the one variable polynomial:
 $$[W_{\bs{\mu},\bs{\sigma}^0}]_e(T):=1+\sum_{j\in \Lambda_e}a_j p^{-i\bs{\tau}_{(n)}^0\cdot  \bs{\alpha}_{(n) j}}T^{\bs{\theta}_{(n)}\cdot  \bs{\alpha}_{(n) j}}=\widetilde{[h]_e}(p^{-i\bs{\tau}_{(n)}^0\cdot  \widehat{\bs{\alpha}}_{(n) e}}T^{\bs{\theta}_{(n)}\cdot  \widehat{\bs{\alpha}}_{(n) e}}).$$

Moreover, if to each root  $c_{\bs{\mu},0}$ of $[W_{\bs{\mu},\bs{\sigma}^0}]_e(T)$ we associate

\begin{equation}\label{afracine-entiere}
c_{\bs{\theta},0}:=c_{\bs{\mu},0}p^{-i\frac{\bs{\tau}_{(n)}^0 \cdot  \bs{\alpha}_{(n) e}}{\bs{\theta}_{(n)}\cdot  \bs{\alpha}_{(n) e}}}.
\end{equation}
Then $\left|c_{\bs{\theta},0}\right|=\left|c_{\bs{\mu},0}\right|$ and $c_{\bs{\theta},0}$ is a root of the polynomial (not depending on $p$ nor on $\bs{\tau}^0$):
\begin{equation}\label{polyhtheta}
 1+\sum_{j\in \Lambda_e}a_j T^{\bs{\theta}_{(n)}\cdot  \bs{\alpha}_{(n) j}}.
\end{equation}

Note that the expression (\ref{polyhtheta}) is really a polynomial according to the property (\ref{entier}) on $\bs{\theta}_{(n)}$.

\begin{rqs}
For $p$ large enough $W_{\bs{\mu},\bs{\sigma}^0}(X,Y)$ is well-defined by putting $X=p^{-1}$ and $Y=p^{-t_{m,\bs{\mu},\bs{\sigma}^0}}$ for $\bs{\tau}_{(n)}^0\in \mathbb{R}^n$  generically chosen by avoiding a countable closed sets of empty interior if necessary.
Indeed  $p^{-1}\in \mathbb{C}\setminus \mathbb{R}_{-}$ and if $p^{-t_{m,\bs{\mu},\bs{\sigma}^0}}\in e^{i\mathfrak{b}}\mathbb{R}_{+}$ we would have the existence of $m_p\in \mathbf{\mathbb{N}^*}$ such that:
\begin{equation}\label{biendef}
% use packages: array
\begin{array}{lll}
\Im(t_{m,\bs{\mu},\bs{\sigma}^0}) & = & \ds\Im\left(-\frac{\log\left(\Omega_{\bs{\mu},\bs{\sigma}^0}(p^{-1})\right)}{\log(p)}+\frac{2 i \pi m}{\log(p)}\right) \\ 
 & = & \ds\frac{-\arg(c_{\bs{\mu},0})+O\left(p^{-\vartheta_{1}}\right)+2 \pi m}{\log(p)} \\ 
 & = & \ds\frac{-\arg(c_{\bs{\theta},0})+O\left(p^{-\vartheta_{1}}\right)+2 \pi m}{\log(p)} - \frac{\bs{\tau}_{(n)}^0 \cdot  \bs{\alpha}_{(n) e}}{\bs{\theta}_{(n)}\cdot  \bs{\alpha}_{(n) e}} = \frac{-\mathfrak{b} + 2 m_p\pi}{\log(p)}. 
\end{array} 
\end{equation}
So, for $m$ fixed, if the equality (\ref{biendef}) was satisfied for infinitely many prime numbers $p$, we could write $\frac{\bs{\tau}_{(n)}^0 \cdot  \bs{\alpha}_{(n) e}}{\bs{\theta}_{(n)}\cdot  \bs{\alpha}_{(n) e}}$ as a limit of a sequence $(S_m)_p$ according to $p$.

But the set $\mathfrak{L}$ of all the limits of these countable sequences is a countable set. Consequently, if we choose $\bs{\tau}_{(n)}^0\in \mathbb{R}^n$ generically so that $\frac{\bs{\tau}_{(n)}^0 \cdot  \bs{\alpha}_{(n) e}}{\bs{\theta}_{(n)}\cdot  \bs{\alpha}_{(n) e}}\notin \mathfrak{L}$ by avoiding a countable union of hypersurfaces of dimension $n-1$ of empty interior, we may assume that $p^{-t_{m,\bs{\mu},\bs{\sigma}^0}}$ for all $m$ and for $p$ large enough.

\end{rqs}

Notice that if $[W_{\bs{\mu},\bs{\sigma}^0}]_e(T)$ is not cyclotomic, then there exists at least one root  $c_{\bs{\mu},0}$ of modulus stricly less than  $1$ which will provide a Puiseux branch  $\Omega_{\bs{\mu},\bs{\sigma}^0}(X)$ satisfying $|\Omega_{\bs{\mu},\bs{\sigma}^0}(X)|<1$ for $|X|$ small enough.

\vspace{0.4cm}

So now let us deal entirely this particular case where  $[W_{\bs{\mu},\bs{\sigma}^0}]_e(T)$ is not cyclotomic: we will prove that there are, among the two factors of  $t \longmapsto Z(\mathbf{s}_{(n)}^0+t\bs{\theta}_{(n)})$ which appear in the writing of Theorem  \ref{afana}, many more zeros coming from the factor $t\longmapsto\prod_{p \leq M_\frac{1}{\nu+1}}W_{\bs{\mu},\bs{\sigma}^0}\left(p^{-1},p^{-t}\right)=\prod_{p \leq M_\frac{1}{\nu+1}}h\left(p^{-s^0_1+t\theta_1},\dots,p^{-s^0_n+t\theta_n},p^{-c}\right)$ (for $\nu>>0$) than poles coming from \\ $t\longmapsto\prod_{\bs{\beta}\in \Nr}\zeta_{M_\frac{1}{\nu+1}} \left((\mathbf{s}^0_{(n)}+t\bs{\theta}_{(n)},c)\cdot\bs{\alpha}\cdot {}^t \! \bs{\beta} \right)^{-\gamma(\bs{\beta})}$ for $t$ lying inside a region $\Delta_{u,\nu,\eta}$ in a neighbourhood on the right of $\Re(t)=0$ determined by (for $\nu,\eta,u>0$):
 
 \begin{center}
 
 \begin{tabular}{ll}
 $\Delta_{u,\nu,\eta}:$ & $\frac{1}{\nu+1}<\Re\left(t\right)< \frac{1}{\nu} $\\ 
  & $0<u<\Im\left(t\right)<u+\eta.$
 \end{tabular}
 \end{center}
 
In this way we will show the accumulation of zeros $t_{m,\bs{\mu},\bs{\sigma}^0}$ ($m\in \mathbb{Z}$, $p$ prime) of $Z\left(\mathbf{s}_{(n)}^0+t\bs{\theta}_{(n)}\right) = Z^{n+1}\left(\left(\mathbf{s}_{(n)}^0,c\right)+t\left(\bs{\theta}_{(n)},0\right)\right)$  inside $\Xi_{u,\eta} = \bigcup_{\nu\geq 1}\Delta_{u,\nu,\eta}$.
 
By the way, we could note that this particular case does not require the use of generic arguments which consist in moving if necessary the parameters  $\bs{\sigma}^0_{(n)}$ or $\bs{\theta}_{(n)}$. Moreover, we can also give an estimation (in function of $\nu$ and $\eta$) of the number of zeros  $t_{m,\bs{\mu},\bs{\sigma}^0}$ inside $\Delta_{u,\nu,\eta}$. It is in this sense that this case is more simple than the case where  $[W_{\bs{\mu},\bs{\sigma}^0}]_e$ is cyclotomic -- the case which will be dealt with later.

 \begin{lemme}\label{afthm_reel}
 
 We suppose that  $[W_{\bs{\mu},\bs{\sigma}^0}]_{e}$ is not cyclotomic.
 
The edge $\partial \mathbf{W}_c(0)$ of $\mathbf{W}_c(0)$ is a natural boundary for  $Z\left(\mathbf{s}\right)=Z^{n+1}(s_1,\dots,s_n,c)$.
 
 In particular, the number $S(\nu,\eta)$ of zeros $t_{m,\bs{\mu},\bs{\sigma}^0}$ of the form   (\ref{afzerotmp})   (counted without their multiplicity) inside the region $\Delta_{\nu,\eta}$ (for $\nu,\eta,u>0$) is such that for all  $N \in \mathbb{N}$:
 
 $$
 S\left(\nu,\eta\right) \geq \frac{\eta \left(C_{0}-1\right)}{\mathcal{K}_N 4\pi}\nu^{N},
 $$
 where $\mathcal{K}_N$ is a constant depending on  $N$ and $C_{0} = \left|{c_{\bs{\mu},0}}^{-1}\right|>1$ is the modulus of the inverse of a root  $c_{\bs{\mu},0}$ of $[W_{\bs{\mu},\bs{\sigma}^0}]_e$ of modulus strictly less than $1$.
  
 \end{lemme}
 
 \begin{dem}
To begin with, let us notice that for $\Re\left(t\right)>\delta$, we have according to (\ref{afana}):
  \begin{displaymath}
 % use packages: array
%  \begin{array}{lll}
 \displaystyle Z\left(\mathbf{s}_{(n)}^0+t\bs{\theta}_{(n)}\right)  = \displaystyle \prod_{p \leq M_\delta}h\left(p^{-s^0_1-t\theta_1},\dots,p^{-s^0_n-t\theta_n},p^{-c}\right)  
    \displaystyle \prod_{\bs{\beta}\in \Nr} \zeta_{M_\delta} \left((\mathbf{s}^0_{(n)}+t\bs{\theta}_{(n)},c)\cdot\bs{\alpha}\cdot {}^t \! \bs{\beta} \right)^{-\gamma(\bs{\beta})}. 
%  \end{array}
 \end{displaymath}
 
 Indeed this writing makes sense because for all $j \in \{1\dots,r\}, (\bs{\sigma}_{(n)}^0+\Re(t)\bs{\theta}_{(n)},c)\cdot  \bs{\alpha}_{\cdot j}\geq \Re(t)\bs{\theta}_{(n)}\cdot  \bs{\alpha}_{(n) j} \geq \delta \ \textrm{according to (\ref{afhyp-theta1})}.
 $

  Consider the zeros and the poles of  $Z\left(\mathbf{s}_{(n)}^0+t\bs{\theta}_{(n)}\right)$ inside the rectangle (for $\nu,\eta,u>0$):
 \begin{center}
 
 \begin{tabular}{ll}
 $\Delta_{u,\nu,\eta}:$ & $\frac{1}{\nu+1}<\Re\left(t\right)< \frac{1}{\nu} $\\ 
  & $0<u<\Im\left(t\right)<u+\eta.$
 \end{tabular}
 \end{center}
 
Firstly let us estimate the number of possible poles inside $\Delta_{\nu,\eta}$ coming from the factor \\ $\ds \prod_{\bs{\beta}\in \Nr} \zeta_{M_{\frac{1}{\nu+1}}} \left((\mathbf{s}^0_{(n)}+t\bs{\theta}_{(n)},c)\cdot\bs{\alpha}\cdot {}^t \! \bs{\beta} \right)^{-\gamma(\bs{\beta})}$.
  Recall (see Theorem \ref{afana}) that $\zeta_{M_{\frac{1}{\nu+1}}}$ has exactly the same zeros and poles as the Riemann zeta function $\zeta$.
 
 If $t_0$ is  such a pole inside $\Delta_{u,\nu,\eta}$, then there exists   $\bs{\beta}\in \Nr$ such that $(\mathbf{s}^0_{(n)}+t_0\bs{\theta}_{(n)},c)\cdot\bs{\alpha}\cdot {}^t \! \bs{\beta}$ is a zero or a pole of the Riemann zeta function $\zeta$; and this quantity satisfies necessarily
 $
 \Re\left(t_0\right) \left((\bs{\theta}_{(n)},0)\cdot \bs{\alpha}\cdot {}^t \! \bs{\beta}\right) \leq \Re\left((\mathbf{s}^0_{(n)}+t_0\bs{\theta}_{(n)},c)\cdot\bs{\alpha}\cdot {}^t \! \bs{\beta}\right) \leq 1.
 $
 
 Consequently we have
 $
 \frac{1}{\nu+1}<\Re\left(t_0\right) \leq \frac{1}{(\bs{\theta}_{(n)},0)\cdot \bs{\alpha}\cdot {}^t \! \bs{\beta}};
 $
 which provides $
 (\bs{\theta}_{(n)},0)\cdot \bs{\alpha}\cdot {}^t \! \bs{\beta} \leq \left(\nu+1\right).
 $
  Furthermore the inequality
 $
 (\bs{\theta}_{(n)},0)\cdot \bs{\alpha}\cdot {}^t \! \bs{\beta} = \sum_{j=1}^{r}\beta_j\left(\bs{\theta}_{(n)}\cdot  \bs{\alpha}_{(n) j}\right) \geq \Vert \beta\Vert \ \left(\textrm{by} \ \left(\ref{afhyp-theta1}\right)\right)
 $ 
 gives:
 \begin{equation}\label{af**}
  \Vert \beta\Vert \leq \left(\nu+1\right).
 \end{equation}
 Moreover $
 \Im\left(t_0\right) < u+\eta
 $
 gives 
 $\Im\left((\mathbf{s}^0_{(n)}+t_0\bs{\theta}_{(n)},c)\cdot\bs{\alpha}\cdot {}^t \! \bs{\beta} \right) 
   =O\left(\left(\nu+1\right)\left(u+\eta\right)\right).$ 
 
 After having fixed $\eta>0$, the number of zeros or poles of one $\zeta$-factor of \\ $ \ds\prod_{\bs{\beta}\in \Nr} \zeta_{M_{\frac{1}{\nu+1}}} \left((\mathbf{s}^0_{(n)}+t\bs{\theta}_{(n)},c)\cdot\bs{\alpha}\cdot {}^t \! \bs{\beta} \right)^{-\gamma(\bs{\beta})}$ is given by $\ds
 O\left(\left(\nu+1\right)\log\left(\nu+1\right)\right),
 $
according to a classical result concerning the estimation of the number of nontrivial zeros of the Riemann zeta function of imaginary part less than $\left(\nu+1\right)$.

 In addition, the same pole can, by $\left(\ref{af**}\right)$, appear in at most  $\left(\nu+1\right)^r$ terms; which gives at most
 $\ds
 O\left(\left(\nu+1\right)^{r+1}\log\left(\nu+1\right)\right)
 $
  poles inside $\Delta_{\nu,\eta}$ (counted without their multiplicity).

\vspace{0.3cm} 

 On the other hand, let us estimate the number of zeros $S\left(\nu,\eta\right)$ coming from  \\ $\prod_{p \leq M_{\frac{1}{\nu+1}}}h\left(p^{-s^0_1-t\theta_1},\dots,p^{-s^0_n-t\theta_n},p^{-c}\right) = \prod_{p\leq M_{\frac{1}{\nu+1}}}W_{\bs{\mu},\bs{\sigma}^0}\left(p^{-1},p^{-t}\right)$ inside $\Delta_{\nu,\eta}$.
 
 We consider for this the Puiseux branches  of $W_{\bs{\mu},\bs{\sigma}^0}(X,Y)$ in the neighbourhood of $X=0$ ($X\in \mathbb{C}\setminus \mathbb{R}_{-}$) of the form (\ref{afbranchepuiseux}).
  We know that the first term $c_{\bs{\mu},0}$ of a branch is a root of  $[W_{\bs{\mu},\bs{\sigma}^0}]_e(T)$ and that reciprocally each root of this polynomial determines the main term of a Puiseux branch (see Proposition \ref{athpuiseuxpolgen} page \pageref{athpuiseuxpolgen}).
 And since by hypothesis $[W_{\bs{\mu},\bs{\sigma}^0}]_e(T)$ is not cyclotomic, then there exists a root  $c_{\bs{\mu},0}$ of modulus stricly less than  $1$.
 
Then consider in particular a Puiseux branch $\Omega_{\bs{\mu},\bs{\sigma}^0}(X)$ having this first term $c_{\bs{\mu},0}$ with $\left|c_{\bs{\mu},0}\right|<1$ and put $C_{0} = |c_{\bs{\mu},0}|^{-1}>1$.
 
 For $p$ prime we write $\ds
 \Omega_{\bs{\mu},\bs{\sigma}^0}\left(p^{-1}\right)= c_{\bs{\mu},0}+c_{\bs{\mu},1}p^{-\vartheta_{1}}+\Omega_{\bs{\mu},\bs{\sigma}^0,2}\left(p^{-1}\right).
 $
 Thus, some zeros of $t\longmapsto W_{\bs{\mu},\bs{\sigma}^0}\left(p^{-1},p^{-t}\right)$ for $p$ prime can be expressed as follows:
 \begin{equation}\label{afformetmp}
 t_{m,\bs{\mu},\bs{\sigma}^0} =  -\frac{\log\left(c_{\bs{\mu},0}+c_{\bs{\mu},1}p^{-\vartheta_{1}}+\Omega_{\bs{\mu},\bs{\sigma}^0,2}\left(p^{-1}\right)\right)}{\log\left(p\right)}+\frac{2\pi m i}{\log\left(p\right)}
 \end{equation}
 where $m \in \mathbb{Z}$.
  To have $t_{m,\bs{\mu},\bs{\sigma}^0} \in \Delta_{\nu,\eta}$, we must have:
 \begin{equation}\label{afreel}
  \frac{1}{\nu+1}<-\frac{\log\left|c_{\bs{\mu},0}+c_{\bs{\mu},1}p^{-\vartheta_{1}}+\Omega_{\bs{\mu},\bs{\sigma}^0,2}\left(p^{-1}\right)\right|}{\log\left(p\right)}<\frac{1}{\nu}.
 \end{equation}
Let us prove that this inequality is well satisfied for $p$ lying in a suitable interval.
 
\begin{equation}\label{areflemme11}
% use packages: array
\begin{array}{l}
\textrm{Firstly we can assume that }   \Re\left(\frac{c_{\bs{\mu},1}}{c_{\bs{\mu},0}}\right)\neq 0:  \textrm{this property which exploits}  \\ 
\textrm{the genericity in} \ \bs{\tau}_{(n)}^0\in \mathbb{R}^n  \textrm{ will be proved in Lemma \ref{afargumentgenerique} page \pageref{afargumentgenerique}.}
\end{array}
\end{equation}

\vspace{0.3cm}
 Thus there exists $p_0 \in \mathbb{N}$ such that for $p>p_0$ we have either:
 \begin{equation}\label{afA}
 \left|1+\frac{c_{\bs{\mu},1}}{c_{\bs{\mu},0}}p^{-\vartheta_{1}}+\frac{\Omega_{\bs{\mu},\bs{\sigma}^0,2}\left(p^{-1}\right)}{c_{\bs{\mu},0}}\right| > 1 \ \ \textrm{if} \ \Re\left(\frac{c_{\bs{\mu},1}}{c_{\bs{\mu},0}}\right) > 0;
 \end{equation}
   or:
 \begin{equation}\label{afB}
  \left|1+\frac{c_{\bs{\mu},1}}{c_{\bs{\mu},0}}p^{-\vartheta_{1}}+\frac{\Omega_{\bs{\mu},\bs{\sigma}^0,2}\left(p^{-1}\right)}{c_{\bs{\mu},0}}\right| < 1 \ \ \textrm{if} \ \Re\left(\frac{c_{\bs{\mu},1}}{c_{\bs{\mu},0}}\right) < 0.
 \end{equation}
 In $\left(\ref{afA}\right)$ or $\left(\ref{afB}\right)$ for $\nu$ large enough and
 \begin{displaymath}
 % use packages: array
 \begin{array}{ll}
 p & >\max\left(\left[4\nu\left|\Re\left(\frac{c_{\bs{\mu},1}}{c_{\bs{\mu},0}}\right)\right|\left(C_{0}^{\frac{\vartheta_{1}}{2}}\right)^{\nu}\right]^{\frac{1}{\vartheta_{1}}},  \left[4\left(\nu+1\right)\left|\Re\left(\frac{c_{\bs{\mu},1}}{c_{\bs{\mu},0}}\right)\right|\left(C_{0}^{\frac{\vartheta_{1}}{2}}\right)^{\nu+1}\right]^{\frac{1}{\vartheta_{1}}}\right) \\ 
  & =\left[4\left(\nu+1\right)\left|\Re\left(\frac{c_{\bs{\mu},1}}{c_{\bs{\mu},0}}\right)\right|\left(C_{0}^{\frac{\vartheta_{1}}{2}}\right)^{\nu+1}\right]^{\frac{1}{\vartheta_{1}}} \\ & (\textrm{which is possible according to Proposition  \ref{athpuiseuxpolgen} page \pageref{athpuiseuxpolgen} because} \\ & |c_{\bs{\mu},0}|=|c_{\bs{\theta},0}|>0 \ \textrm{and} \  |c_{\bs{\mu},1}| \ \textrm{is bounded independently of} \ p), \\
  & \textrm{we obtain}: 
 \end{array}
 \end{displaymath}
  \begin{equation}\label{afineg}
  \left(-1\right)^{\varepsilon}\left|1+\frac{c_{\bs{\mu},1}}{c_{\bs{\mu},0}}p^{-\vartheta_{1}}+\frac{\Omega_{\bs{\mu},\bs{\sigma}^0,2}\left(p^{-1}\right)}{c_{\bs{\mu},0}}\right|^{-\nu-\varepsilon} < \frac{1+\left(-1\right)^{\varepsilon}\left(C_{0}^{\frac{\vartheta_{1}}{2}}\right)^{\nu+\varepsilon}}{\left(C_{0}^{\frac{\vartheta_{1}}{2}}\right)^{\nu+\varepsilon}}
 \end{equation}
 where $\varepsilon \in \left\{0,1\right\}$.

  Indeed, the inequality $\left(\ref{afA}\right)$ gives:
 \begin{enumerate}
  \item for $\varepsilon=0:$
 $\ds
 \left|1+\frac{c_{\bs{\mu},1}}{c_{\bs{\mu},0}}p^{-\vartheta_{1}}+\frac{\Omega_{\bs{\mu},\bs{\sigma}^0,2}\left(p^{-1}\right)}{c_{\bs{\mu},0}}\right|^{-\nu} < 1 < \frac{1+\left(C_{0}^{\frac{\vartheta_{1}}{2}}\right)^{\nu}}{\left(C_{0}^{\frac{\vartheta_{1}}{2}}\right)^{\nu}};
 $
 \item for $\varepsilon=1,$ since $p>\left[4\left(\nu+1\right)\Re\left(\frac{c_{\bs{\mu},1}}{c_{\bs{\mu},0}}\right)\left(C_{0}^{\frac{\vartheta_{1}}{2}}\right)^{\nu+1}\right]^{\frac{1}{\vartheta_{1}}}$:
 
 \begin{displaymath}
 % use packages: array
 \begin{array}{ll}
 & \log \left|1+\frac{c_{\bs{\mu},1}}{c_{\bs{\mu},0}}p^{-\vartheta_{1}}+\frac{\Omega_{\bs{\mu},\bs{\sigma}^0,2}\left(p^{-1}\right)}{c_{\bs{\mu},0}}\right|^{-\nu-1}  \\
 =   &\left(\left(1+\frac{c_{\bs{\mu},1}}{c_{\bs{\mu},0}}p^{-\vartheta_{1}}+\frac{\Omega_{\bs{\mu},\bs{\sigma}^0,2}\left(p^{-1}\right)}{c_{\bs{\mu},0}}\right)\overline{\left(1+\frac{c_{\bs{\mu},1}}{c_{\bs{\mu},0}}p^{-\vartheta_{1}}+\frac{\Omega_{\bs{\mu},\bs{\sigma}^0,2}\left(p^{-1}\right)}{c_{\bs{\mu},0}}\right)}\right)^{\frac{-\nu-1}{2}}     \\ 
 = & \log \left(1+ 2\Re\left(\frac{c_{\bs{\mu},1}}{c_{\bs{\mu},0}} \right)p^{-\vartheta_{1}}+o\left(p^{-\vartheta_{1}}\right) \right)^{\frac{-\nu-1}{2}}  \\ 
 = & \frac{-\nu-1}{2}\left(2\Re\left(\frac{c_{\bs{\mu},1}}{c_{\bs{\mu},0}} \right)p^{-\vartheta_{1}}+o\left(p^{-\vartheta_{1}}\right) \right) \\ 
 > &  -2\left(\nu+1\right)\Re\left(\frac{c_{\bs{\mu},1}}{c_{\bs{\mu},0}} \right)p^{-\vartheta_{1}} \ \textrm{for} \ \nu \ \textrm{large enough} \ (\nu\geq \nu_0)\\
 > & -\frac{1}{2}\frac{1}{C_{0}^{\frac{\vartheta_{1}}{2}\left(\nu+1\right)}} \\
 > & \log\left(1- \frac{1}{C_{0}^{\frac{\vartheta_{1}}{2}\left(\nu+1\right)}}\right) = -\frac{1}{C_{0}^{\frac{\vartheta_{1}}{2}\left(\nu+1\right)}} + o\left(\frac{1}{C_{0}^{\frac{\vartheta_{1}}{2}\left(\nu+1\right)}}\right);
 \end{array}
 \end{displaymath}
which provides the inequality desired (for $\nu \geq \nu_0$):

  $\ds
 \left|1+\frac{c_{\bs{\mu},1}}{c_{\bs{\mu},0}}p^{-\vartheta_{1}}+\frac{\Omega_{\bs{\mu},\bs{\sigma}^0,2}\left(p^{-1}\right)}{c_{\bs{\mu},0}}\right|^{-\nu-1} > 1-\left(C_{0}^{\frac{\vartheta_{1}}{2}}\right)^{-1-\nu} = \frac{-1+\left(C_{0}^{\frac{\vartheta_{1}}{2}}\right)^{1+\nu}}{\left(C_{0}^{\frac{\vartheta_{1}}{2}}\right)^{1+\nu}}.
 $
 \end{enumerate}
 
 Similarly, the inequality $\left(\ref{afB}\right)$ gives:
 
 \begin{enumerate}
  \item for $\varepsilon=0$, since $p>\left[-4\nu\Re\left(\frac{c_{\bs{\mu},1}}{c_{\bs{\mu},0}}\right)\left(C_{0}^{\frac{\vartheta_{1}}{2}}\right)^{\nu}\right]^{\frac{1}{\vartheta_{1}}}$:
 
 \begin{displaymath}
 % use packages: array
 \begin{array}{lll}
   \log \left|1+\frac{c_{\bs{\mu},1}}{c_{\bs{\mu},0}}p^{-\vartheta_{1}}+\frac{\Omega_{\bs{\mu},\bs{\sigma}^0,2}\left(p^{-1}\right)}{c_{\bs{\mu},0}}\right|^{-\nu}  & = & \frac{\nu}{2}\left(-2\Re\left(\frac{c_{\bs{\mu},1}}{c_{\bs{\mu},0}}\right)p^{-\vartheta_{1}} + o\left(p^{-\vartheta_{1}} \right) \right) \\ 
 &< & 2\nu\left(-\Re\left(\frac{c_{\bs{\mu},1}}{c_{\bs{\mu},0}}\right)p^{-\vartheta_{1}} \right) \ \textrm{for} \ \nu \ \textrm{large enough} \ (\nu \geq \nu_1) \\ 
 &< & \frac{1}{2C_{0}^{\frac{\vartheta_{1}}{2}\nu}}  
 < \log\left(1+\frac{1}{C_{0}^{\frac{\vartheta_{1}}{2}\nu}} \right) = \frac{1}{C_{0}^{\frac{\vartheta_{1}}{2}\nu}} + o\left(\frac{1}{C_{0}^{\frac{\vartheta_{1}}{2}\nu}} \right);
 \end{array}
 \end{displaymath}
 which guarantees (for $\nu\geq \nu_1$) $
 \left|1+\frac{c_{\bs{\mu},1}}{c_{\bs{\mu},0}}p^{-\vartheta_{1}}+\frac{\Omega_{\bs{\mu},\bs{\sigma}^0,2}\left(p^{-1}\right)}{c_{\bs{\mu},0}}\right|^{-\nu} < 1+\frac{1}{C_{0}^{\frac{\vartheta_{1}}{2}\nu}}.
 $
 
 \item for $\varepsilon=1:$
 $\ds
 \left|1+\frac{c_{\bs{\mu},1}}{c_{\bs{\mu},0}}p^{-\vartheta_{1}}+\frac{\Omega_{\bs{\mu},\bs{\sigma}^0,2}\left(p^{-1}\right)}{c_{\bs{\mu},0}}\right|^{-\nu-1} > 1 > 1-\frac{1}{C_{0}^{\frac{\vartheta_{1}}{2}(\nu+1)}} = \frac{C_{0}^{\frac{\vartheta_{1}}{2}(\nu+1)}-1}{C_{0}^{\frac{\vartheta_{1}}{2}(\nu+1)}}.
 $
 \end{enumerate}

 Now if we choose $\ds C_{0}^{\nu}\left(1+\frac{1}{C_{0}^{\nu\frac{\vartheta_{1}}{2}}}\right) \leq p \leq C_{0}^{\nu+1}\left(1-\frac{1}{C_{0}^{\left(\nu+1\right)\frac{\vartheta_{1}}{2}}}\right),$ (which is compatible with the previous condition on $p$ to have $\left(\ref{afineg}\right)$) then $\left(\ref{afreel}\right)$ occurs since according to $\left(\ref{afineg}\right)$ we have:
  \begin{center}
  \begin{tabular}{lllll}
  $C_{0}^{\nu}\left|1+\frac{c_{\bs{\mu},1}}{c_{\bs{\mu},0}}p^{-\vartheta_{1}}+\frac{\Omega_{\bs{\mu},\bs{\sigma}^0,2}\left(p^{-1}\right)}{c_{\bs{\mu},0}}\right|^{-\nu}$ & $< C_{0}^{\nu}\left(1+\frac{1}{C_{0}^{\nu\frac{\vartheta_{1}}{2}}}\right)$ \\
  & $ \leq p $ \\
  & $\leq C_{0}^{\nu+1}\left(1-\frac{1}{C_{0}^{\left(\nu+1\right)\frac{\vartheta_{1}}{2}}}\right)$ \\ 
  & $< C_{0}^{\nu+1}\left|1+\frac{c_{\bs{\mu},1}}{c_{\bs{\mu},0}}p^{-\vartheta_{1}}+\frac{\Omega_{\bs{\mu},\bs{\sigma}^0,2}\left(p^{-1}\right)}{c_{\bs{\mu},0}}\right|^{-\nu-1};$ \\ 
 \end{tabular}
 \end{center}
 
and finally by taking the logarithm of both sides we deduce $\left(\ref{afreel}\right)$.

 Now, $\eta>0$ being fixed, if we choose $\nu$ as a positive integer such that  $\frac{2\pi}{\log\left(C_{0}^{\nu}+1\right)}<\eta$, then for all prime number $p$ such that
 $\ds
 C_{0}^{\nu}\left(1+\frac{1}{C_{0}^{\nu\frac{\vartheta_{1}}{2}}}\right) \leq p \leq C_{0}^{\nu+1}\left(1-\frac{1}{C_{0}^{\left(\nu+1\right)\frac{\vartheta_{1}}{2}}}\right),
 $ 
 we will have $t_{m,\bs{\mu},\bs{\sigma}^0} \in \Delta_{\nu,\eta}$ if and only if $\ds
 u<\frac{2\pi m}{\log\left(p\right)} - \frac{\arg\left(\Omega_{\bs{\mu},\bs{\sigma}^0}\left(p^{-1}\right)\right)}{\log\left(p\right)} < u+\eta,
 $
 which is equivalent to:
 \begin{equation}\label{afboite2}
 \frac{u\log\left(p\right)}{2\pi} + \frac{\arg\left(\Omega_{\bs{\mu},\bs{\sigma}^0}\left(p^{-1}\right)\right)}{2\pi}<m<\frac{\left(u+\eta\right)\log\left(p\right)}{2\pi} + \frac{\arg\left(\Omega_{\bs{\mu},\bs{\sigma}^0}\left(p^{-1}\right)\right)}{2\pi}.
 \end{equation}
Hence, for a fixed  $p$, we will have   $\frac{\eta\log\left(p\right)}{2\pi}+\varpi$ zeros $t_{m,\bs{\mu},\bs{\sigma}^0}$ of $W_{\bs{\mu},\bs{\sigma}^0}\left(p^{-1},p^{-t}\right)$ inside $\Delta_{\nu,\eta}$ where $|\varpi|\leq 1$.

 Finally, if $S^{*}\left(\nu,\eta\right)$ denotes the number of zeros of   $\prod_{p \leq M_{\frac{1}{\nu+1}}}h\left(p^{-s^0_1-t\theta_1},\dots,p^{-s^0_n-t\theta_n},p^{-c}\right)$ inside $\Delta_{\nu,\eta}$ a priori counted with their multiplicity, we will have:
 \begin{equation}\label{afS}
 S^{*}\left(\nu,\eta\right) \geq \sum_{C_{0}^{\nu}\left(1+\frac{1}{C_{0}^{\nu\frac{\vartheta_{1}}{2}}}\right) \leq p \leq C_{0}^{\nu+1}\left(1-\frac{1}{C_{0}^{\left(\nu+1\right)\frac{\vartheta_{1}}{2}}}\right)}\left(\frac{\eta\log\left(p\right)}{2\pi}+\varpi\right).
 \end{equation}
 
 By taking $\nu$ large enough so that $C_{0}^{-\nu\frac{\vartheta_{1}}{2}} < \frac{C_{0}-1}{2\left(C_{0}^{1-\frac{\vartheta_{1}}{2}}+1\right)}$ and by using the prime number theorem  (i.e. $\displaystyle{\sum_{p \leq x}\log\left(p\right) \sim x}$), the estimation  $\left(\ref{afS}\right)$ gives:
 
 \begin{eqnarray*}
  S^{*}\left(\nu,\eta\right) &\geq& \frac{ C_{0}^{\nu}\eta\left(C_{0}-1\right)}{4\pi}-\sum_{C_{0}^{\nu}\left(1+\frac{1}{C_{0}^{\nu\frac{\vartheta_{1}}{2}}}\right) \leq p \leq C_{0}^{\nu+1}\left(1-\frac{1}{C_{0}^{\left(\nu+1\right)\frac{\vartheta_{1}}{2}}}\right)}1  \\
   &\geq& \frac{ C_{0}^{\nu}\eta\left(C_{0}-1\right)}{4\pi} - \frac{C_{0}^{\nu+1}}{\log\left(C_{0}^{\nu+1}\right)}  \\
   &\sim& \frac{ C_{0}^{\nu}\eta\left(C_{0}-1\right)}{4\pi}.
 \end{eqnarray*}
 
For being able to minorate  $S\left(\nu,\eta\right)$, we want to majorate the multiplicity of a zero or a pole  $t_{m,\bs{\mu},\bs{\sigma}^0}$.
 Thus given a prime number $p$ and an integer  $m$, we want to majorate:
  $$\mathcal{M}\left(m,p\right)=\#\left\{\left(m',p'\right) \mid m' \in \mathbb{Z}, p' \ \textrm{prime,} \ t_{m,\bs{\mu},\bs{\sigma}^0}=t_{m',\bs{\mu}',\bs{\sigma}^{0}} \ \textrm{with} \ \bs{\mu'}=(p',\bs{\tau}^0_{(n)},\bs{\theta_{(n)}})\right\}.$$
 
Notice that we can suppose without loss of generality that if $p'$ is such that there exists an integer $m$ such that $t_{m,\bs{\mu},\bs{\sigma}^0}=t_{m',\bs{\mu}',\bs{\sigma}^0}$, then $p' \geq p$.
 
 In addition we have:
 \begin{equation}\label{aflogg}
 \begin{array}{c}
  -\log \Omega_{\bs{\mu},\bs{\sigma}^0}\left(p^{-1}\right) = -\log\left(c_{\bs{\mu},0}\right) + O\left(p^{-\vartheta_{1}}\right); \\
 -\log \Omega_{\bs{\mu}',\bs{\sigma}^0}\left(p'^{-1}\right) = -\log\left(c_{\bs{\mu}',0}\right) + O\left(p^{-\vartheta_{1}}\right).
 \end{array}
 \end{equation}
 Furthermore by (\ref{afracine-entiere}) we can write, for all prime number $p$,
 $\ds
 c_{\bs{\mu},0} = c_{\bs{\theta},0}p^{i\frac{\bs{\tau}_{(n)}^0 \cdot  \bs{\alpha}_{(n) e}}{\bs{\theta}_{(n)}\cdot  \bs{\alpha}_{(n) e}}},
 $
 where $c_{\bs{\theta},0}$ does not depend on $p$.
  We remark also that $\Re\left(\log\left(c_{\bs{\theta},0}\right)\right)=\log|c_{\bs{\theta},0}| \neq 0$ because $|c_{\bs{\theta},0}|=\left|c_{\bs{\mu},0}\right|<1$.
  According to $\left(\ref{aflogg}\right)$, the equality $t_{m,\bs{\mu},\bs{\sigma}^0}=t_{m',\bs{\mu}',\bs{\sigma}^0}$ provides:
 
 \begin{equation}\label{afcaca}
 \frac{-\log\left(c_{\bs{\theta},0}\right) + O\left(p^{-\vartheta_{1}}\right)}{\log\left(p\right)}+\frac{2 i \pi m}{\log\left(p\right)} = \frac{-\log\left(c_{\bs{\theta},0}\right) + O\left(p^{-\vartheta_{1}}\right)}{\log\left(p'\right)}+\frac{2 i \pi m'}{\log\left(p'\right)}.
 \end{equation}
 
By identifying the real and the imaginary parts of $\left(\ref{afcaca}\right)$, we obtain the estimations:
 \begin{equation*}
 \begin{cases}
 \displaystyle -\log|c_{\bs{\theta},0}|\left(\frac{1}{\log\left(p\right)}-\frac{1}{\log\left(p'\right)}\right) & \displaystyle =  O\left(\frac{1}{p^{\vartheta_{1}}\log\left(p\right)}\right), \\
 \displaystyle -\arg\left(c_{\bs{\theta},0}\right)\left(\frac{1}{\log\left(p\right)}-\frac{1}{\log\left(p'\right)}\right) + 2\pi\left(\frac{m}{\log\left(p\right)}-\frac{m'}{\log\left(p'\right)}\right) & \displaystyle =  O\left(\frac{1}{p^{\vartheta_{1}}\log\left(p\right)}\right).
 \end{cases}
 \end{equation*}
 
 And since $\log|c_{\bs{\theta},0}| \neq 0$, we have:
 \begin{equation}\label{afqq}
  \begin{cases}
  \displaystyle \frac{1}{\log\left(p\right)}-\frac{1}{\log\left(p'\right)} & = \displaystyle O\left(\frac{1}{p^{\vartheta_{1}}\log\left(p\right)}\right), \\
  \displaystyle \frac{m}{\log\left(p\right)}-\frac{m'}{\log\left(p'\right)} & = \displaystyle O\left(\frac{1}{p^{\vartheta_{1}}\log\left(p\right)}\right).
  \end{cases}
 \end{equation}
 
 The first line of $\left(\ref{afqq}\right)$ permits to claim that $\log\left(p'\right)-\log\left(p\right) = O\left(\frac{\log\left(p'\right)}{p^{\vartheta_{1}}}\right).$
 Consequently there exists an absolute constant $A_1$ such that if $p'$ is such that there exists $m'$ verifying $t_{m',\bs{\mu}',\bs{\sigma}^0}=t_{m,\bs{\mu},\bs{\sigma}^0}$ then $
 \log\left(p'\right)-\log\left(p\right) \leq A_1\frac{\log\left(p'\right)}{p^{\vartheta_{1}}}.
 $ 
 So we have $\ds
 \log\left(p'\right) \leq \frac{\log\left(p\right)}{1-\frac{A_1}{p^{\vartheta_{1}}}} \leq \log\left(p\right)\left(1 + \frac{A_2}{p^{\vartheta_{1}}}\right);
 $
 where $A_2$ is an absolute constant (we can for example choose $A_2=2A_1$).
  If there exists $m'$ such that $t_{m',\bs{\mu}',\bs{\sigma}^0}=t_{m,\bs{\mu},\bs{\sigma}^0}$, then $p'$ satisfies necessarily
 \begin{equation}\label{afencadr}
  p' \leq p^{1+\frac{A_2}{p^{\vartheta_{1}}}}.
 \end{equation}
 
 For a fixed $p$, let us count the number   $\mathcal{M'}\left(p\right)$ of $p'$ satisfying $\left(\ref{afencadr}\right)$.
 For this we use the prime number theorem which gives the following estimation $
 \pi\left(x\right) = \int_{2}^{x}\frac{d t}{\log\left(t\right)} + O\left(x e^{-\mathfrak{m}\sqrt{\log\left(x\right)}}\right)
 $ (where $\mathfrak{m}$ is an explicit absolute constant) for the number of prime numbers $\pi\left(x\right)$ less than  $x$.
  Hence we obtain:

 $\ds \mathcal{M'}\left(p\right) = \displaystyle\pi\left(p^{1+ A_2 p^{-\vartheta_{1}}}\right)-\pi\left(p\right)  
   = \displaystyle\int_{2}^{p^{1+ A_2 p^{-\vartheta_{1}}}}\frac{d t}{\log\left(t\right)} + O\left(p e^{-\mathfrak{m}\sqrt{\log\left(p\right)}}\right).$
 
 But we have uniformly in  $t \in \left[p, p^{1+ A_2 p^{-\vartheta_{1}}} \right]$ $\ds
 \log\left(t\right) = \log\left(p\right) + O\left(\log\left(p\right) p^{-\vartheta_{1}}\right);
 $
 which provides:
 \begin{displaymath}
 % use packages: array
 \begin{array}{lll}
 \mathcal{M'}\left(p\right) &=& \displaystyle\frac{1}{\log\left(p\right) + O\left(\log\left(p\right) p^{-\vartheta_{1}}\right)}\left(p^{1+ A_2 p^{-\vartheta_{1}}}-p\right) + O\left(p e^{-\mathfrak{m}\sqrt{\log\left(p\right)}}\right) \\ 
  &=& \displaystyle O\left(\frac{p}{\log\left(p\right)}\left(p^{A_2 p^{-\vartheta_{1}}}-1\right)\right) +  O\left(p e^{-\mathfrak{m}\sqrt{\log\left(p\right)}}\right)\\ 
  &=& \displaystyle O\left(\frac{p}{\log\left(p\right)}\left(e^{A_2 \log\left(p\right) p^{-\vartheta_{1}}}-1\right)\right) + O\left(p e^{-\mathfrak{m}\sqrt{\log\left(p\right)}}\right) \\ 
  &=& \displaystyle O\left(p^{1-\vartheta_{1}}\right) +  O\left(p e^{-\mathfrak{m}\sqrt{\log\left(p\right)}}\right)\\ 
  &=& \displaystyle O\left(p e^{-\mathfrak{m}\sqrt{\log\left(p\right)}}\right).
 \end{array}
 \end{displaymath}
 
 Now, having fixed an integer  $m \in \mathbb{Z}$ and a prime number $p$, let us consider a prime number $p'$ verifying $\left(\ref{afencadr}\right)$  and let us estimate the number of integers  $m'$ such that $t_{m,\bs{\mu},\bs{\sigma}^0} = t_{m',\bs{\mu}',\bs{\sigma}^0}$.
 
According to $\left(\ref{afqq}\right)$, we have
 $
 \displaystyle \frac{m}{\log\left(p\right)}-\frac{m'}{\log\left(p'\right)}  = \displaystyle O\left(\frac{1}{p^{\vartheta_{1}}\log\left(p\right)}\right).
 $ 
 But since $p'$ verifies $\left(\ref{afencadr}\right)$, we have $\ds
 \log\left(p'\right) = \log\left(p\right) + O\left(\log\left(p\right) p^{-\vartheta_{1}} \right);
 $ 
 and consequently:
  \begin{displaymath}
 % use packages: array
 \begin{array}{lll}
 \displaystyle m-m'\frac{\log\left(p\right)}{\log\left(p'\right)} &=&\displaystyle O\left(p^{-\vartheta_{1}}\right) \\ 
 \displaystyle m-m'\left(\frac{1}{1+O\left(p^{-\vartheta_{1}}\right)}\right) &=&\displaystyle O\left(p^{-\vartheta_{1}}\right) \\ 
 \displaystyle m-m'\left(1+O\left(p^{-\vartheta_{1}}\right)\right) &=&\displaystyle O\left(p^{-\vartheta_{1}}\right) \\ 
 \displaystyle m-m' &=&\displaystyle O\left(p^{-\vartheta_{1}}\right) + O\left(m' p^{-\vartheta_{1}}\right).
 \end{array}
 \end{displaymath}
 
 Moreover, if $t_{m',\bs{\mu}',\bs{\sigma}^0} \in \Delta_{\nu,\eta}$, then by $\left(\ref{afboite2}\right)$ $m'$ must verify $\ds
 m' = O\left(\log\left(p'\right)\right) = O\left(\log\left(p\right)\right);
 $
 and hence $\ds
 m-m' = O\left(\log\left(p\right) p^{-\vartheta_{1}}\right).
 $ 
 In particular, for $p$ large enough, $p>p_1$ ($p_1$ being an absolute constant), we deduce $\ds
 |m-m'|<\frac{1}{2};
 $
 and then
 $
 m=m'.
 $
 Hence if  $p>p_1$, the couples $\left(m',p'\right)$ such that $t_{m',\bs{\mu}',\bs{\sigma}^0} = t_{m,\bs{\mu},\bs{\sigma}^0}$ are necessarily such that  $m=m'$.
  And finally
 $\ds
 \mathcal{M}\left(m,p\right) = \mathcal{M'}\left(p\right) = O\left(p e^{-\mathfrak{m}\sqrt{\log\left(p\right)}}\right).
 $ 

 As a conclusion, if $p$ is such that $C_{0}^{\nu}\left(1+\frac{1}{C_{0}^{\nu\frac{\vartheta_{1}}{2}}}\right) \leq p \leq C_{0}^{\nu+1}\left(1-\frac{1}{C_{0}^{\left(\nu+1\right)\frac{\vartheta_{1}}{2}}}\right)$, then for all $N \in \mathbb{N}$, there exists in particular a constant $\mathcal{K}_N$ which depends on $N$ such that for all $m \in \mathbb{N}$,
 $\ds
 \mathcal{M}\left(m,p\right) \leq \mathcal{K}_N \frac{C_{0}^{\nu}}{\nu^{N}}.
 $ 
 Thus for all  $N \in \mathbb{N}$, we have finally $$
 S\left(\nu,\eta\right) \geq \frac{S^{*}\left(\nu,\eta\right)}{\mathcal{K}_N \frac{C_{0}^{\nu}}{\nu^{N}}} \sim \frac{\eta \left(C_{0}-1\right)}{\mathcal{K}_N 4\pi}\nu^{N}.
 $$ 
 For $N>r+1$, we have in particular $\left(\nu+1\right)^{r+1}\log\left(\nu+1\right) = o\left(S\left(\nu,\eta\right)\right)$ when $\nu$ tends to infinity; which completes the proof of this lemma.  \CQFD

 \end{dem}
 
 \vspace{0.5cm}

 \textit{The case where $[W_{\bs{\mu},\bs{\sigma}^0}]_e(T)$ is not cyclotomic being now completely dealt with, we assume from now on that   $[W_{\bs{\mu},\bs{\sigma}^0}]_e(T)$ is cyclotomic.}

The problem is more complicated when the polynomial $[W_{\bs{\mu},\bs{\sigma}^0}]_e(T)$ is cyclotomic.

In this case, it is necessary to consider the second term of the Puiseux branches to prove that there exists at least one satisfying the desired condition $|\Omega_{\bs{\mu},\bs{\sigma}^0}(X)|<1$. Moreover, to have this condition, it is also necessary to refine our choice of direction $\bs{\theta}$. 
Indeed, by choosing a suitable  $\bs{\theta}$, we will see that we can find two Puiseux series of opposite initial term $\pm c_{\bs{\mu},0}$ with the same second term $c_{\bs{\mu},1}X^{\vartheta_{1}}$.
In this way, although it is not possible to have  $\left|c_{\bs{\mu},0}\right|<1$ since here $[W_{\bs{\mu}, \bs{\sigma}^0}]_e$ is cyclotomic, one of these two branches will be of modulus stricly less than  $1$ for $X>0$ small enough whenever   $\arg\left(\frac{c_{\bs{\mu},1}}{c_{\bs{\mu},0}}\right) \neq \frac{\pi}{2} \mod(\pi)$. And we use an argument of genericity in the imaginary part  $\bs{\tau}^0 = \left(\bs{\tau}_{(n)}^0,0\right)$ of $\mathbf{s}^0$ to ensure 
$\ds
\arg\left(\frac{c_{\bs{\mu},1}}{c_{\bs{\mu},0}}\right) \neq \frac{\pi}{2} \mod(\pi).
$

\vspace{0.3cm}

We will use later the following lemma to justify the existence of a particular index $e'\in \{1,\dots,r\}$ which will be crucial in the computation of the second term of the Puiseux branches.

\begin{lemme}\label{aW-noncyclo}
Consider an half-line $e^{i\mathfrak{b}}\mathbb{R}_{+}$ and a corresponding determination of the logarithm so that  $W_{\bs{\mu},\bs{\sigma}^0}\mid_{\bs{\tau}^0=\textbf{0}}(X,Y)$ is defined for  $X\in \mathbb{C}\setminus \mathbb{R}_{-}$ and $Y\in \mathbb{C}\setminus e^{i\mathfrak{b}}\mathbb{R}_{+}$.
Suppose that there exists $\mathfrak{c}\in \mathbb{C}\setminus e^{i\mathfrak{b}}\mathbb{R}_{+}$ such that for all $X\in \mathbb{C}\setminus \mathbb{R}_{-}$ and for $\mathbf{s}_{(n)}^0$  generically chosen in $\mathcal{B}\cap \partial \mathbf{W}_c(0)$ we have
$
W_{\bs{\mu},\bs{\sigma}^0}\mid_{\bs{\tau}^0=\textbf{0}}\left(X,\mathfrak{c}\right) = 0.
$
Then necessarily $\left|\mathfrak{c}\right|\neq 1$.
\end{lemme}

\begin{dem}
 Suppose that there exists $\mathfrak{c}\in \mathbb{C}\setminus e^{i\mathfrak{b}}\mathbb{R}_{+}$ such that for all $X\in \mathbb{C}\setminus \mathbb{R}_{-}:$ 
$$\ds
W_{\bs{\mu},\bs{\sigma}^0}\mid_{\bs{\tau}^0=\textbf{0}}(X,\mathfrak{c})=0.
$$
Then we will have:
\begin{equation}\label{ahtnul}
% use packages: array
\begin{array}{lll}
\displaystyle W_{\bs{\mu},\bs{\sigma}^0}\mid_{\bs{\tau}^0=\textbf{0}}\left(X,\mathfrak{c}\right) & = & \displaystyle 1+\sum_{j=1}^{r}a_j X^{\bs{\sigma}^0\cdot  \bs{\alpha}_{\cdot j}}\mathfrak{c}^{\bs{\theta}\cdot  \bs{\alpha}_{\cdot j}} \\ 
 & = & \displaystyle h\left(X^{\sigma^0_1}\mathfrak{c}^{\theta_1},\dots, X^{\sigma^0_n}\mathfrak{c}^{\theta_n}\right) = 0. 
\end{array}
\end{equation}

Moreover, the generic choice of  $\bs{\sigma}^0_{(n)}$ (see (\ref{afdcond})) means that the only constraint its components must verify is  $\bs{\sigma}^0\cdot  \bs{\alpha}_{\cdot e} = 0$.

Consequently, since we have supposed without loss of generality that $\alpha_{n, e}\neq 0$ (see page \pageref{afhyp-theta1}),  we can consider $\bs{\sigma}^0_{(n)}\in \mathbb{R}^n$ as a $(n-1)$-tuple $\widetilde{\bs{\sigma}}^0=(\widetilde{\sigma}^0_1,\dots,\widetilde{\sigma}^0_{n-1})\in U\subseteq \mathbb{R}^{n-1}$ ($U$ being an open set of $\mathbb{R}^{n-1}$) by putting:

\begin{equation*}
 \begin{cases}
 \displaystyle \sigma^0_{\ell} = \widetilde{\sigma}^0_{\ell} & (\ell \in \{1,\dots,n-1\}), \\
 \displaystyle \sigma^0_n = -\frac{1}{\alpha_{n, e}}\left(\sum_{i=1}^{n-1}\alpha_{i, e}\widetilde{\sigma}^0_i + c  \ \alpha_{n+1, e}\right) & .
 \end{cases}
 \end{equation*}

Then for all $x\in \mathbb{R}_{>0}$ $(x\neq 1)$ let us define:

\begin{displaymath}
% use packages: array
\begin{array}{cccc}
 & U & \longrightarrow & \mathbb{R}^{n} \\ 
\Phi_{x}: & & & \\
 & \widetilde{\bs{\sigma}}^0=(\widetilde{\sigma}^0_1,\dots,\widetilde{\sigma}^0_{n-1}) & \longmapsto & \left(x^{\widetilde{\sigma}^0_1},\dots,x^{\widetilde{\sigma}^0_{n-1}},x^{c}\right).
\end{array}
\end{displaymath}

It is clear that $\bigcup_{x>0}\Phi_{x}(U)$ describes a nonempty open set $(0,\infty)^{n}$.

Moreover, there exists a nonempty open set $U'\subseteq \mathbb{R}^{n-1}\times \left(\mathbb{R}_{>0}\setminus \{1\}\right)$ such that for all $(y_1,\dots, y_n)\in U'$ there exists $(\widetilde{\sigma}^0_1,\dots,\widetilde{\sigma}^0_{n-1})\in U$ and $x>0$ ($x\neq 1$) such that  $$(y_1^{\alpha_{n, e}},\dots,y_n^{\alpha_{n, e}})=\Phi_{\exp\left(\frac{\alpha_{n, e}}{c}\log(y_n)\right)}\left(c\frac{\log(y_1)}{\log(y_n)},\dots,c\frac{\log(y_{n-1})}{\log(y_n)}\right) \in \textrm{Im}\left(\Phi_{\exp\left(\frac{\alpha_{n, e}}{c}\log(y_n)\right)}\right).$$

Thus for all $(y_1,\dots, y_n)\in U'$ we have:

\begin{displaymath}
% use packages: array
\begin{array}{lll}
\displaystyle h\left(y_1^{\alpha_{n, e}} \mathfrak{c}^{\theta_1},\dots,y_{n-1}^{\alpha_{n, e}} \mathfrak{c}^{\theta_{n-1}}, \mathfrak{c}^{\theta_{n}}y_n^{-\alpha_{n+1, e}}\prod_{\ell=1}^{n-1}y_{\ell}^{-\alpha_{\ell, e}},y_n^{\alpha_{n, e}}\right) & = & \displaystyle h\left(x^{\sigma^0_1}\mathfrak{c}^{\theta_1},\dots, x^{\sigma^0_n}\mathfrak{c}^{\theta_n},x^c\right) \\
& & \textrm{with} \ x=\exp\left(\frac{\alpha_{n, e}}{c}\log(y_n)\right) \\ 
 & = & 0 \ \textrm{according to (\ref{ahtnul})}.
\end{array}
\end{displaymath}

But the function 
$\ds
\left(y_1,\dots,y_{n}\right) \longmapsto \displaystyle h\left(y_1^{\alpha_{n, e}} \mathfrak{c}^{\theta_1},\dots,y_{n-1}^{\alpha_{n, e}} \mathfrak{c}^{\theta_{n-1}}, \mathfrak{c}^{\theta_{n}}y_n^{-\alpha_{n+1, e}}\prod_{\ell=1}^{n-1}y_{\ell}^{-\alpha_{\ell, e}},y_n^{\alpha_{n, e}}\right)
$
is holomorphic on $(\mathbb{C}^{*})^{n}$.
And since it vanishes on an open set $U'$ of $(0,\infty)^{n}$, we have in fact
$$
\textrm{for all} \ (y_1,\dots,y_{n})\in (\mathbb{C}^{*})^{n}, \  h\left(y_1^{\alpha_{n, e}} \mathfrak{c}^{\theta_1},\dots,y_{n-1}^{\alpha_{n, e}} \mathfrak{c}^{\theta_{n-1}}, \mathfrak{c}^{\theta_{n}}y_n^{-\alpha_{n+1, e}}\prod_{\ell=1}^{n-1}y_{\ell}^{-\alpha_{\ell, e}},y_n^{\alpha_{n, e}}\right)=0.
$$
Hence the polynomial $ h(X_1,\dots,X_{n+1})$ vanishes on $\mathcal{H}\cap (\mathbb{C}^{*})^{n}$ where $\mathcal{H}$ is the complex hypersurface defined by the equation $\ds
\mathfrak{c}^{-\bs{\theta}_{(n)}\cdot  \bs{\alpha}_{(n) e}}\textbf{X}^{\bs{\alpha}_{\cdot e}}-1 = \prod_{\ell=1}^{n+1}X_{\ell}^{\alpha_{\ell, e}} \mathfrak{c}^{-\theta_{\ell}\alpha_{\ell, e}} -1 =0.
$
We deduce that the polynomial $X_1 \cdots X_{n}h(X_1,\dots,X_{n+1})$ vanishes on the whole hypersurface $\mathcal{H}$ and hence the polynomial $\mathfrak{c}^{-\bs{\theta}_{(n)}\cdot  \bs{\alpha}_{(n) e}} \textbf{X}^{\bs{\alpha}_{\cdot e}} -1$
divides a power of the polynomial $X_1\dots X_{n} ~h(X_1,\dots,X_{n+1})$. 
Since the polynomials $\mathfrak{c}^{-\bs{\theta}_{(n)}\cdot  \bs{\alpha}_{(n) e}} \textbf{X}^{\bs{\alpha}_{\cdot e}} -1$ and $X_1\dots X_{n}$ are relatively prime, we deduce that the polynomial  
\begin{equation}\label{adiviseurh}
P_{\mathfrak{c}}(\textbf{X}):= \mathfrak{c}^{-\bs{\theta}_{(n)}\cdot  \bs{\alpha}_{(n) e}}\textbf{X}^{\bs{\alpha}_{\cdot e}}-1
\end{equation}
 necessarily divides a power of $h$; and hence $P_{\mathfrak{c}}(\textbf{X})$ divides also $h$ because all irreducible factors of $P_{\mathfrak{c}}(\textbf{X})$ are of multiplicity $1$.

And since $h$ is with rational coefficients and  $\mathfrak{c}$ is an algebraic number, the polynomial
$$
Q(\textbf{X}):=\prod_{\mathfrak{c}'}P_{\mathfrak{c}'}(\textbf{X}) \in \mathbb{Q}[\textbf{X}]
$$
(where the product is done over all the conjugates $\mathfrak{c}'$ of $\mathfrak{c}$)
also divides $h$.

Remark that $Q(\textbf{X})$ can be reduced in fact to an one variable polynomial (by the change of variable $T:=\textbf{X}^{\bs{\alpha}_{\cdot e}}$).

Now, if we assume by absurd that $|\mathfrak{c}|=1$, then the polynomial $Q(\textbf{X})$, having rational coefficients, is necessarily cyclotomic because all its roots would be of modulus $1$; which is not possible since  by hypothesis $h$ does not contain any cyclotomic factor.
\CQFD
\end{dem}

Now consider a particular index $e'\in \{1,\dots,r\}$ verifying: 
\begin{equation}\label{afchoixe}
  \bs{\sigma}^0\cdot \bs{\alpha}_{\cdot e'}>0 \ \textrm{is minimal}
\end{equation}
among the indexes $j_0$ such that $\bs{\sigma}^0\cdot \bs{\alpha}_{\cdot j_0}>0$ and which verify:
 $$\sum_{\{j:\bs{\alpha}_{\cdot j}-\bs{\alpha}_{\cdot j_0}\in \mathbb{Q}\bs{\alpha}_{\cdot e}\}}a_j p^{-i\bs{\tau}^0_{(n)}\cdot \bs{\alpha}_{(n) j}}{c_{\bs{\mu},0}}^{\bs{\theta}_{(n)}\cdot  \bs{\alpha}_{(n) j}}\neq 0.$$

The existence of such index $e'$ is ensured by the previous lemma.

Indeed, consider a set  $J$ made of representatives of each class of the following equivalence relation  $\sim$:
$$
\bs{\alpha}_{\cdot j} \sim \bs{\alpha}_{\cdot j'} \Longleftrightarrow \bs{\alpha}_{\cdot j}-\bs{\alpha}_{\cdot j'}\in \mathbb{Q}\bs{\alpha}_{\cdot e}.
$$

And write

$$
W_{\bs{\mu},\bs{\sigma}^0}(X,Y) = [W_{\bs{\mu},\bs{\sigma}^0}]_e(Y) + \sum_{j_0\in J; j_0\not\sim e}X^{\bs{\sigma}^0\cdot \bs{\alpha}_{\cdot j_0}}R_{\bs{\mu}, j_0}(Y);
$$

where
$$
R_{\bs{\mu}, j_0}(Y)=\sum_{j\sim j_0}a_j p^{-i\bs{\tau}^0_{(n)}\cdot \bs{\alpha}_{(n) j}} Y^{\bs{\theta}_{(n)}\cdot  \bs{\alpha}_{(n) j}}.
$$

Recall that since we suppose here that  $[W_{\bs{\mu},\bs{\sigma}^0}]_e(Y)$ is a cyclotomic polynomial, all its roots are of modulus  $1$.

Now let $c_{\bs{\mu},0}$ be a root of $[W_{\bs{\mu},\bs{\sigma}^0}]_e$.

We have:
$$
R_{\bs{\mu}, j_0}(c_{\bs{\mu},0}) = \sum_{j \sim j_0} a_j p^{-i\bs{\tau}_{(n)}^0\cdot \left(\bs{\alpha}_{(n) j}-\bs{\alpha}_{(n) e}\frac{\bs{\theta}_{(n)}\cdot  \bs{\alpha}_{(n) j}}{\bs{\theta}_{(n)}\cdot  \bs{\alpha}_{(n) e}}\right)} {c_{\bs{\theta},0}}^{\bs{\theta}_{(n)}\cdot  \bs{\alpha}_{(n) j}} \ \textrm{according to}\ (\ref{afracine-entiere}).
$$

But if $\bs{\alpha}_{\cdot j} = \bs{\alpha}_{\cdot j_0}+q \bs{\alpha}_{\cdot e}$, we obtain:

\begin{displaymath}
% use packages: array
\begin{array}{lll}
\bs{\alpha}_{(n) j}-\bs{\alpha}_{(n) e}\frac{\bs{\theta}_{(n)}\cdot  \bs{\alpha}_{(n) j}}{\bs{\theta}_{(n)}\cdot  \bs{\alpha}_{(n) e}} & = & \bs{\alpha}_{(n) j_0} + q \bs{\alpha}_{(n) e} -\bs{\alpha}_{(n) e}\frac{\bs{\theta}_{(n)}\cdot\bs{\alpha}_{(n) j_0}+q\bs{\theta}_{(n)}\cdot  \bs{\alpha}_{(n) e}}{\bs{\theta}_{(n)}\cdot  \bs{\alpha}_{(n) e}} \\ 
 & = & \bs{\alpha}_{(n) j_0}-\bs{\alpha}_{(n) e}\frac{\bs{\theta}_{(n)}\cdot\bs{\alpha}_{(n) j_0}}{\bs{\theta}_{(n)}\cdot  \bs{\alpha}_{(n) e}}.
\end{array}
\end{displaymath}

Consequently since the $\bs{\alpha}_{(n) j}-\bs{\alpha}_{(n) e}\frac{\bs{\theta}_{(n)}\cdot  \bs{\alpha}_{(n) j}}{\bs{\theta}_{(n)}\cdot  \bs{\alpha}_{(n) e}}$ are all equal for  $j\sim j_0$  we have $\ds
R_{\bs{\mu}, j_0}(c_{\bs{\mu},0}) = 0
$
is equivalent to: 
\begin{equation}\label{aequationnulle?}
R_{j_0}(c_{\bs{\theta},0}):=\sum_{j\sim j_0} a_j  {c_{\bs{\theta},0}}^{\bs{\theta}_{(n)}\cdot  \bs{\alpha}_{(n) j}} = 0.
\end{equation}

Thus if $e'$ does not exist, we would have for all $X\in \mathbb{C}\setminus \mathbb{R}_{-}$:
$$
W_{\bs{\mu},\bs{\sigma}^0}\mid_{\bs{\tau}^0=\textbf{0}}(X,c_{\bs{\theta},0}) = [h_{\theta}]_e(c_{\bs{\theta},0}) + \sum_{j_0\in J; j_0\not\sim e} X^{\bs{\sigma}^0\cdot \bs{\alpha}_{\cdot j_0}}R_{j_0}(c_{\bs{\theta},0})=0;
$$

which is impossible according to Lemma \ref{aW-noncyclo} since here $|c_{\bs{\theta},0}|=|c_{\bs{\mu},0}|=1$.

Obviously, it is possible to have some $j_0$ such that
\begin{equation}\label{aminimalindice}
\bs{\sigma}^0\cdot \bs{\alpha}_{\cdot e'} = \bs{\sigma}^0\cdot \bs{\alpha}_{\cdot j_0}.
\end{equation}

However, if $\bs{\sigma}^0\in \mathbb{R}^n$ is  generically chosen so that ${\mathbf{s}^0}\in \mathcal{B}\cap \partial \mathbf{W}(0)$, the equality (\ref{aminimalindice}) implies necessarily that $j_0\sim e'$. 

\vspace{0.3cm}

 Throughout what follows, we need the direction $\bs{\theta} = \left(\bs{\theta}_{(n)},0\right)\in \mathbb{Q}^n\times \{0\}$ to verify, in addition of (\ref{afhyp-theta1}), the two following conditions:

\begin{equation}\label{afhyp-theta2}
% use packages: array
\begin{array}{l}
  \bs{\theta}_{(n)}\cdot\bs{\widehat{\alpha}}_{(n) e} \in \mathbb{Z}_{+} \ \textrm{is even};   \\
  \bs{\theta}_{(n)}\cdot\bs{\alpha}_{(n) e'} \in \mathbb{Z}_{+} \ \textrm{is odd}.
\end{array}
\end{equation}

Note that, \textit{although the vectors  $\bs{\alpha}_{\cdot e}$ and $\bs{\alpha}_{\cdot e'}$ are not collinear  (according to (\ref{afdcond}) since $\bs{\sigma}^0\cdot  \bs{\alpha}_{\cdot e}=0$ and $\bs{\sigma}^0\cdot \bs{\alpha}_{\cdot e'}>0$), it would be possible to have  $\bs{\alpha}_{(n) e'}\in \mathbb{Q}\bs{\alpha}_{(n) e}$}; and in this case it would not be possible to choose such  $\bs{\theta}$ with $\theta_{n+1}=0$ and verifying (\ref{afhyp-theta2}).

To escape this difficulty,  \textit{we use the hypothesis (H)} of Theorem \ref{afprincipal} which, because we know that  $\bs{\alpha}_{\cdot e}$ and $\bs{\alpha}_{\cdot e'}$ are not collinear, ensures that $\bs{\alpha}_{(n) e'}\notin \mathbb{Q}\bs{\alpha}_{(n) e}$; and consequently it is possible to find  $\bs{\theta} = \left(\bs{\theta}_{(n)},0\right)$ verifying (\ref{afhyp-theta2}).

\vspace{0.3cm}

The two following lemmas are crucial to prove the existence of an accumulation of zeros $t_{m,\bs{\mu},\bs{\sigma}^0}\in \Xi_{u,\eta}$ of positive real part.

\begin{lemme}\label{afargumentgenerique}
 Let $\Omega_{\bs{\mu},\bs{\sigma}^0}(X) = c_{\bs{\mu},0} + c_{\bs{\mu},1}X^{\vartheta_{1}} + o\left(X^{\vartheta_{1}}\right)$ be a Puiseux branch of initial term $c_{\bs{\mu},0}$, root of $[W_{\bs{\mu},\bs{\sigma}^0}]_e$.
Moving generically $\bs{\tau}_{(n)}^0\in\mathbb{R}^n$ so that $\mathbf{s}_{(n)}^0\in \mathcal{B}\cap \partial \mathbf{W}_c(0)$ if necessary, we can assume 
$$
\arg\left(\frac{c_{\bs{\mu},1}}{c_{\bs{\mu},0}}\right) \neq \frac{\pi}{2} \mod(\pi).
$$
\end{lemme}

\begin{rqs}
 This lemma does not require the hypothesis that $[W_{\bs{\mu},\bs{\sigma}^0}]_e(T)$ is cyclotomic; moreover this result is used in the proof of Lemma \ref{afthm_reel} page \pageref{afthm_reel} (see (\ref{areflemme11}) page \pageref{areflemme11}).
\end{rqs}

This lemma will be proved just after the following lemma:

\begin{lemme}\label{aaccumulation-zeros-h}
Assume that $[W_{\bs{\mu},\bs{\sigma}^0}]_e$ is a cyclotomic polynomial.

  There exists a Puiseux series $\Omega_{\bs{\mu},\bs{\sigma}^0}(X)$, solution of $W_{\bs{\mu},\bs{\sigma}^0}(X,Y)=0$ and verifying
$$
|\Omega_{\bs{\mu},\bs{\sigma}^0}(X)|<1 \ \textrm{for} \ X>0 \ \textrm{small enough};
$$
which provides an infinite number of zeros $t_{m,\bs{\mu},\bs{\sigma}^0} \in \Xi_{u,\eta}$ of $\prod_{p\leq M_{\delta}}h(p^{-s^0_1-t\theta_1},\dots, p^{-s^0_n-t\theta_n},p^{-c})$ as $\delta$ tends to $0$. 

\end{lemme}

\begin{dem}
Consider a Puiseux branch that we will write: 
$$\Omega^{+}_{\bs{\mu}, \bs{\sigma}^0}\left(X\right) = c_{\bs{\mu},0} + c^{+}_{\bs{\mu}, 1}X^{\vartheta^{+}_{1}}+o\left(X^{\vartheta^{+}_{1}}\right);$$
 of main term the root  $c_{\bs{\mu},0}$ of $[W_{\bs{\mu},\bs{\sigma}^0}]_e$ introduced previously which is of multiplicity $m_0= 1$ since $\mathcal{F}(\bs{\alpha}_{\cdot e})$ is non-degenerate in the sense of Definition \ref{anouvcond2}.

According to Lemma \ref{afargumentgenerique}, moving $\bs{\tau}_{(n)}^0$ if necessary, we can suppose that
$
\arg\left(\frac{c^{+}_{\bs{\mu}, 1}}{c_{\bs{\mu},0}}\right) \neq \frac{\pi}{2} \mod(\pi).
$
Thus we have  
$
\frac{\pi}{2} < \arg\left(\frac{c^{+}_{\bs{\mu}, 1}}{c_{\bs{\mu},0}}\right) < \frac{3 \pi}{2} \ \textrm{or} \ \frac{\pi}{2} < \arg\left(-\frac{
c^{+}_{\bs{\mu}, 1}}{c_{\bs{\mu},0}}\right) < \frac{3 \pi}{2}.
$
Since $[W_{\bs{\mu},\bs{\sigma}^0}]_e$ is supposed to be cyclotomic, the main term of  $\Omega^{+}_{\bs{\mu}, \bs{\sigma}^0}$ is of modulus $|c_{\bs{\mu},0}|=1$.

But if we assume firstly that \ $\frac{\pi}{2} < \arg\left(\frac{c^{+}_{\bs{\mu}, 1}}{c_{\bs{\mu},0}}\right) < \frac{3 \pi}{2}$, we have:

$$\left|\Omega^{+}_{\bs{\mu}, \bs{\sigma}^0}(X)\right| = \left|1 + \frac {c^{+}_{\bs{\mu}, 1}}{c_{\bs{\mu},0}} X^{\vartheta^{+}_{1}} + o(X^{\vartheta^{+}_{1}})\right| \ < \ 1 \ \ \  \text{{\it  for $X>0$ small.}}$$
Hence $\Omega_{\bs{\mu},\bs{\sigma}^0} (X) = \Omega^{+}_{\bs{\mu}, \bs{\sigma}^0} (X)$ suits and it is the Puiseux series that we have looked for.

Now assume that
$\ds
\frac{\pi}{2} < \arg\left(-\frac{c^{+}_{\bs{\mu}, 1}}{c_{\bs{\mu},0}}\right) < \frac{3 \pi}{2}.
$

We will show that the particular choice of  $\bs{\theta}$ (see (\ref{afhyp-theta2}) page \pageref{afhyp-theta2}) permits to find a Puiseux series with initial term  $-c_{\bs{\mu},0}$ and with the same second term as that of $\Omega^{+}_{\bs{\mu}, \bs{\sigma}^0}$ and which hence will be the series we have looked for.

Firstly, since $\bs{\theta}_{(n)}\in \mathbb{Q}^n$ has been chosen so that, for  $j \in \Lambda_e$, $\bs{\theta}\cdot  \bs{\alpha}_{\cdot j}$ is even, we know that  $-c_{\bs{\mu},0}$ is also a root of  $[W_{\bs{\mu},\bs{\sigma}^0}]_e$; furthermore notice that this root $-c_{\bs{\mu},0}$ provides a corresponding Puiseux series solution of  $W_{\bs{\mu}, \bs{\sigma}^0}(X,Y)=0$ according to Proposition \ref{athpuiseuxpolgen}.

So consider the following Puiseux branch that we will write:
$$\Omega^{-}_{\bs{\mu}, \bs{\sigma}^0}\left(X\right) = -c_{\bs{\mu},0} + c^{-}_{\bs{\mu}, 1}X^{\vartheta^{-}_{1}}+o\left(X^{\vartheta^{-}_{1}}\right);$$
and whose initial term is this root $-c_{\bs{\mu},0}$.
Let us compare the two terms $c^{-}_{\bs{\mu}, 1}X^{\vartheta^{-}_{1}}$ and $c^{+}_{\bs{\mu}, 1}X^{\vartheta^{+}_{1}}$.
 We use for this the fact that the terms of lowest degree in $X$ of $W_{\bs{\mu},\bs{\sigma}^0}\left(X,\Omega^{\pm}_{\bs{\mu}, \bs{\sigma}^0}(X)\right)$ cancel each other; and these terms coincide with those of $W_{\bs{\mu},\bs{\sigma}^0}\left(X,\pm c_{\bs{\mu},0} + c^{\pm}_{\bs{\mu}, 1}X^{\vartheta^{\pm}_{1}}\right)$. 
And these terms of lowest degree are also those of the following expression:
\begin{equation}\label{aequasimpl}
c^{\pm}_{\bs{\mu}, 1} [W_{\bs{\mu},\bs{\sigma}^0}]'_e\left(\pm c_{\bs{\mu},0}\right) X^{\vartheta^{\pm}_{1}} +  X^{\bs{\sigma}^0\cdot \bs{\alpha}_{\cdot e'}} R_{\bs{\mu}, e'}\left(\pm c_{\bs{\mu},0}\right).
\end{equation}                 

Then on one hand concerning the branch $\Omega^{+}_{\bs{\mu}, \bs{\sigma}^0}$ we have:
\begin{equation}\label{aalgopuiseux1}
 c^{+}_{\bs{\mu}, 1} = -\frac{R_{\bs{\mu}, e'}(c_{\bs{\mu},0})}{[W_{\bs{\mu},\bs{\sigma}^0}]'_e(c_{\bs{\mu},0})};
\end{equation}
and on the other hand concerning the branch $\Omega^{-}_{\bs{\mu}, \bs{\sigma}^0}$ we have:
\begin{equation}\label{aalgopuiseux2}
 c^{-}_{\bs{\mu}, 1} = -\frac{R_{\bs{\mu}, e'}(-c_{\bs{\mu},0})}{[W_{\bs{\mu},\bs{\sigma}^0}]'_e(-c_{\bs{\mu},0})} = -\frac{(-1)^{\bs{\theta}_{(n)}\cdot\bs{\alpha}_{(n) e'}}R_{\bs{\mu}, e'}(c_{\bs{\mu},0})}{-[W_{\bs{\mu},\bs{\sigma}^0}]'_e(c_{\bs{\mu},0})}.
\end{equation}

And since $\bs{\theta}_{(n)}\cdot\bs{\alpha}_{(n) e'}$ is an odd integer we obtain  $c^{+}_{\bs{\mu}, 1} = c^{-}_{\bs{\mu}, 1}$.

Hence finally there exists a Puiseux series:  
$$\Omega_{\bs{\mu},\bs{\sigma}^0}(X) = -c_{\bs{\mu},0} + c^{+}_{\bs{\mu}, 1}X^{\vartheta^{+}_{1}} + o\left(X^{\vartheta^{+}_{1}}\right);$$
which is such that $|\Omega_{\bs{\mu},\bs{\sigma}^0}(X)|<1$ for $X>0$ small enough.
This series provides some zeros: 
$$
t_{m,\bs{\mu},\bs{\sigma}^0} = -\frac{\log\left(\Omega_{\bs{\mu},\bs{\sigma}^0}\left(p^{-1}\right)\right)}{\log(p)} + \frac{2 \pi m i}{\log(p)};
$$
where $m\in \mathbb{Z}$ and $p$ is a prime number. And we will have $t_{m,\bs{\mu},\bs{\sigma}^0} \in \Xi_{u,\eta}$ if $u<\Im(t_{m,\bs{\mu},\bs{\sigma}^0})<u+\eta$; meaning that if: 

$$
u<\frac{2\pi m}{\log\left(p\right)} - \frac{\arg\left(\Omega_{\bs{\mu}, \bs{\sigma}^0}\left(p^{-1}\right)\right)}{\log\left(p\right)} < u+\eta,
$$
which is equivalent to:
\begin{equation}\label{aboite}
\frac{u\log\left(p\right)}{2\pi} + \frac{\arg\left(\Omega_{\bs{\mu}, \bs{\sigma}^0}\left(p^{-1}\right)\right)}{2\pi}<m<\frac{\left(u+\eta\right)\log\left(p\right)}{2\pi} + \frac{\arg\left(\Omega_{\bs{\mu}, \bs{\sigma}^0}\left(p^{-1}\right)\right)}{2\pi}.
\end{equation}
Hence we will have for $p$ large enough some zeros of  $t \to W_{\bs{\mu}, \bs{\sigma}^0}(p^{-1}, p^{-t})$ inside $\Xi_{u,\eta}$.
And there exists infinitely many zeros  $t_{m,\bs{\mu},\bs{\sigma}^0}\in \Xi_{u,\eta}$ of $\ds \prod_{p\leq M_{\delta}}h(p^{-s^0_1-t\theta_1},\dots,p^{-s^0_n-t\theta_n},p^{-c})$ when $\delta$ tends to $0$; which completes the proof of this lemma.
\CQFD
\end{dem}

Now let us prove Lemma \ref{afargumentgenerique}:

\begin{dem}[of Lemma \ref{afargumentgenerique}]
To begin with, we identify the dependence on $p$ and $\bs{\tau}_{(n)}^0$ of $c_{\bs{\mu},1}$.

According to Lemma \ref{aaccumulation-zeros-h}  we have:
$$
c_{\bs{\mu},1} = \frac{-\sum_{j\mid \bs{\alpha}_{\cdot j}-\bs{\alpha}_{\cdot e'}\in \mathbb{Q}\bs{\alpha}_{\cdot e}} a_j p^{-i\bs{\tau}^0_{(n)}\cdot \bs{\alpha}_{(n) j}}\left(c_{\bs{\mu},0}\right)^{\bs{\theta}_{(n)}\cdot\bs{\alpha}_{(n) j}}}{[W_{\bs{\mu},\bs{\sigma}^0}]'_e\left(c_{\bs{\mu},0}\right)}.
$$ 
Recall that  by (\ref{afracine-entiere}) we have $\ds
c_{\bs{\mu},0} = c_{\bs{\theta},0} p^{i\frac{\bs{\tau}_{(n)}^0 \cdot  \bs{\alpha}_{(n) e}}{\bs{\theta}_{(n)}\cdot  \bs{\alpha}_{(n) e}}};
$
which permits to  identify clearly the dependence on $p$ and $\bs{\tau}_{(n)}^0$ of $c_{\bs{\mu},0}$.
Let us observe in particular the denominator $c_{\bs{\mu},0} [W_{\bs{\mu},\bs{\sigma}^0}]'_e\left(c_{\bs{\mu},0}\right)$ of $\frac{c_{\bs{\mu},1}}{c_{\bs{\mu},0}}$ and let us prove that it does not depend neither on $\bs{\tau}_{(n)}^0$, nor on $p$.
Indeed:

\begin{footnotesize}\begin{displaymath}
% use packages: array
\begin{array}{lll}
\displaystyle c_{\bs{\mu},0} [W_{\bs{\mu},\bs{\sigma}^0}]'_e\left(c_{\bs{\mu},0}\right) & = & \displaystyle c_{\bs{\theta},0} p^{i\frac{\bs{\tau}_{(n)}^0 \cdot  \bs{\alpha}_{(n) e}}{\bs{\theta}_{(n)}\cdot  \bs{\alpha}_{(n) e}}} \\
 & & \displaystyle\sum_{j\in \Lambda_e} a_j (\bs{\theta}\cdot  \bs{\alpha}_{\cdot j}) c_{\bs{\theta},0}^{\bs{\theta}_{(n)}\cdot  \bs{\alpha}_{(n) j}-1}p^{i\left(\frac{\bs{\tau}_{(n)}^0 \cdot  \bs{\alpha}_{(n) e}}{\bs{\theta}_{(n)}\cdot  \bs{\alpha}_{(n) e}}\left(\bs{\theta}_{(n)}\cdot  \bs{\alpha}_{(n) j} -1\right)-\bs{\tau}^0_{(n)}\cdot \bs{\alpha}_{(n) j}\right)} \\ 
 & = & \displaystyle \sum_{j\in \Lambda_e}a_j (\bs{\theta}\cdot  \bs{\alpha}_{\cdot j}) c_{\bs{\theta},0}^{\bs{\theta}_{(n)}\cdot  \bs{\alpha}_{(n) j}},
\end{array}
\end{displaymath}                 \end{footnotesize}

since $j \in \Lambda_e$ implies that  $\frac{\bs{\tau}_{(n)}^0 \cdot  \bs{\alpha}_{(n) e}}{\bs{\theta}_{(n)}\cdot  \bs{\alpha}_{(n) e}}(\bs{\theta}_{(n)}\cdot  \bs{\alpha}_{(n) j}) = \bs{\tau}^0_{(n)}\cdot \bs{\alpha}_{(n) j}.$

Now assume by absurd that for all $\bs{\tau}_{(n)}^0$ in some open ball so that $\mathbf{s}_{(n)}^0\in \mathcal{B}\cap \partial \mathbf{W}_c(0)$ and for all prime number $p$
$\arg\left(\frac{c_{\bs{\mu},1}}{c_{\bs{\mu},0}}\right) = \frac{\pi}{2} \mod(\pi).$
Then we would have
\begin{equation*}
   \arg\left(\frac{c_{\bs{\mu},1}}{c_{\bs{\mu},0}}\right)\in \left\{ \frac{\pi}{2};\frac{3\pi}{2} \right\}.
 \end{equation*}
Put for $k$ such that $\bs{\alpha}_{\cdot k} \in \bs{\alpha}_{\cdot e'} + \mathbb{Q}\bs{\alpha}_{\cdot e}$:
$$
\lambda_{k, \bs{\theta}, \bs{\beta}, \rho}:= \frac{- \ a_k \ c_{\bs{\theta},0}^{\bs{\theta}_{(n)}\cdot\bs{\alpha}_{(n) k}}}{\sum_{j\in \Lambda_e}a_j  (\bs{\theta}\cdot  \bs{\alpha}_{\cdot j}) c_{\bs{\theta},0}^{\bs{\theta}_{(n)}\cdot\bs{\alpha}_{(n) j}}} \in \mathbb{C}
$$
not depending neither on $p$ nor on $\bs{\tau}_{(n)}^0$.
Then we have:

\begin{displaymath}
% use packages: array
\begin{array}{lll}
\displaystyle \frac{c_{\bs{\mu},1}}{c_{\bs{\mu},0}} & = & \displaystyle \sum_{\{k :  \bs{\alpha}_{\cdot k}-\bs{\alpha}_{\cdot e'}\in \mathbb{Q}\bs{\alpha}_{\cdot e}\}}\lambda_{k, \bs{\theta}, \bs{\beta}, \rho} p^{i\left(\bs{\tau}_{(n)}^0\cdot\bs{\alpha}_{(n) e}\frac{\bs{\theta}_{(n)}\cdot\bs{\alpha}_{(n) k}}{\bs{\theta}_{(n)}\cdot  \bs{\alpha}_{(n) e}}-\bs{\tau}_{(n)}^0\cdot\bs{\alpha}_{(n) k}\right)} \\ 
 & = &\displaystyle \sum_{\{k :  \bs{\alpha}_{\cdot k}-\bs{\alpha}_{\cdot e'}\in \mathbb{Q}\bs{\alpha}_{\cdot e}\}}\lambda_{k, \bs{\theta}, \bs{\beta}, \rho} p^{i\bs{\tau}_{(n)}^0\cdot \mathbf{w}_k},
\end{array}
\end{displaymath}

if we write
$
\mathbf{w}_{k}:=\frac{\bs{\theta}_{(n)}\cdot\bs{\alpha}_{(n) k}}{\bs{\theta}_{(n)}\cdot  \bs{\alpha}_{(n) e}}\bs{\alpha}_{(n) e} - \bs{\alpha}_{(n) k}.
$

Remark that these $\mathbf{w}_{k}$ are all equal.
 Indeed if $k$ and $k'$ are such that $\bs{\alpha}_{\cdot k}-\bs{\alpha}_{\cdot e'}\in \mathbb{Q}\bs{\alpha}_{\cdot e}$ and $\bs{\alpha}_{\cdot k'}-\bs{\alpha}_{\cdot e'}\in \mathbb{Q}\bs{\alpha}_{\cdot e}$, then $\bs{\alpha}_{\cdot k}-\bs{\alpha}_{\cdot k'}\in \mathbb{Q}\bs{\alpha}_{\cdot e}$ and hence there exists $q\in \mathbb{Q}$ such that
$
\bs{\alpha}_{\cdot k}-\bs{\alpha}_{\cdot k'} = q  \bs{\alpha}_{\cdot e}.
$
In particular we also have $\bs{\alpha}_{(n) k}-\bs{\alpha}_{(n) k'} = q  \bs{\alpha}_{(n) e}$; consequently
$
\frac{\bs{\theta}_{(n)}\cdot(\bs{\alpha}_{(n) k}-\bs{\alpha}_{(n) k'})}{\bs{\theta}_{(n)}\cdot  \bs{\alpha}_{(n) e}}\bs{\alpha}_{(n) e} = q  \bs{\alpha}_{(n) e} = \bs{\alpha}_{(n) k}-\bs{\alpha}_{(n) k'}
$
and hence $\mathbf{w}_{k}=\mathbf{w}_{k'}=\mathbf{w}_{e'}$.

\vspace{0.3cm}

 \textit{Let us prove that  $\mathbf{w}_{e'}\neq \textbf{0}:$ this property is crucial if we want to move   $\arg\left(\frac{c_{\bs{\mu},1}}{c_{\bs{\mu},0}}\right)$ by choosing $\bs{\tau}_{(n)}^0$ generically. It is exactly here that we use again the hypothesis (H) of theorem  \ref{afprincipal}.}

\vspace{0.3cm}

Indeed, $\mathbf{w}_{e'}=\textbf{0}$ is equivalent to $\bs{\alpha}_{(n) e'} \in \mathbb{Q}\bs{\alpha}_{(n) e}$. But we know that the vectors  $\bs{\alpha}_{\cdot e}$ and $\bs{\alpha}_{\cdot e'}$ are not collinear since $\bs{\sigma}^0\cdot \bs{\alpha}_{\cdot e'}>0$ and $\bs{\sigma}^0\cdot  \bs{\alpha}_{\cdot j}=0$ if and only if $\bs{\alpha}_{\cdot j} \in \mathbb{Q}\bs{\alpha}_{\cdot e}$ according to the property (\ref{afhyp00}) on $\bs{\sigma}^0=\left(\bs{\sigma}^0_{(n)},c\right)$.
And consequently the hypothesis (H) ensures $\bs{\alpha}_{(n) e'} \notin \mathbb{Q}\bs{\alpha}_{(n) e}$ and hence $\mathbf{w}_{e'}\neq \textbf{0}$.

Finally we have
$\ds
\frac{c_{\bs{\mu},1}}{c_{\bs{\mu},0}} = p^{i\bs{\tau}_{(n)}^0\cdot\mathbf{w}_{e'}}\sum_{k\mid \bs{\alpha}_{\cdot k}-\bs{\alpha}_{\cdot e'}\in \mathbb{Q}\bs{\alpha}_{\cdot e}}\lambda_{k, \bs{\theta}, \bs{\beta}, \rho}.
$
Now if we put:
$$
\varphi := \arg\left(\sum_{k\mid \bs{\alpha}_{\cdot k}-\bs{\alpha}_{\cdot e'}\in \mathbb{Q}\bs{\alpha}_{\cdot e}}\lambda_{k, \bs{\theta}, \bs{\beta}, \rho}\right);
$$
($\varphi$ not depending neither on $p$ nor on $\bs{\tau}_{(n)}^0$), then we obtain:

\begin{equation*}
  \bs{\tau}_{(n)}^0\in M:=\bigcup_{p}\bigcup_{m\in \mathbb{Z}}\left\{(\bs{\tau}_{(n)}^0\cdot\mathbf{w}_{e'})\log(p) + \varphi =\frac{\pi}{2} + m\pi  \right\}.
 \end{equation*}

But $M$ is a countable union of affine hyperplanes in $\bs{\tau}_{(n)}^0$ which are of empty interior inside  $\mathbb{R}^n$ (\textit{because $\mathbf{w}_{e'}\neq \textbf{0}$}); and according to Baire's theorem the countable union of these hyperplanes is also of empty interior inside  $\mathbb{R}^n$.
Consequently these previous conditions cannot be satisfied for all  $\bs{\tau}_{(n)}^0$ inside an open ball of  $\mathbb{R}^n$; and we obtain a contradiction to the hypothesis above, which completes the proof of this lemma.   \CQFD

\end{dem}

\vspace{0.3cm}
At this stage we have proved the existence of infinitely many zeros $t_{m,\bs{\mu},\bs{\sigma}^0}\in \Xi_{u,\eta}$ of $t\longmapsto \prod_{p\leq M_{\delta}}h\left(p^{-s^0_1-t\theta_1},\dots, p^{-s^0_n-t\theta_n},p^{-c}\right)$ as $\delta$ tends to $0$.
To prove Theorem \ref{afprincipal} it remains to verify that the accumulation of these zeros $t_{m,\bs{\mu},\bs{\sigma}^0}$ lying inside $\Xi_{u,\eta}$ is not cancelled by possible poles coming from the $\zeta$-factors of  $\prod_{\bs{\beta}\in \Nr} \zeta_{M_\delta} \left((\mathbf{s}^0_{(n)}+t\bs{\theta}_{(n)},c)\cdot\bs{\alpha}\cdot {}^t \! \bs{\beta} \right)^{-\gamma(\bs{\beta})} $ which appear in the writing of Theorem \ref{afana}.

To start with, let us remark that, since $\zeta_{M_\delta}$ has exactly the same zeros and pole as the Riemann zeta function, the possible poles which could cancel the previous zeros $t_{m,\bs{\mu},\bs{\sigma}^0}$ are of the form:

$$
t\left(\beta,\rho\right) = \frac{\rho-\sum_{j=1}^{r}{\beta_j(\mathbf{s}^0\cdot\bs{\alpha}_{\cdot j})}}{\sum_{j=1}^{r}\beta_j (\bs{\theta}\cdot  \bs{\alpha}_{\cdot j})},
$$
where $\bs{\beta}\in \Nr$ and $\rho$ is a zero or a pole of $\zeta$.

In what follows we will prove that, by moving  $\mathbf{s}_{(n)}^0\in \mathcal{B}\cap \partial \mathbf{W}_c(0)$ if necessary, there is at most a finite number of such  $t(\bs{\beta},\rho)$ inside the region   $\Xi_{u,\eta}$ which vanish the factors $t\longmapsto h\left(p^{-s^0_1-t\theta_1},\dots, p^{-s^0_n-t\theta_n},p^{-c}\right)$ ($p$ being any prime number); and consequently they cannot cancel the accumulation of $t_{m,\bs{\mu},\bs{\sigma}^0}\in \Xi_{u,\eta}$.

Then  for all $\bs{\beta}\in \Nr$, for all prime number $p$ and for all $\rho$ zero or pole of $\zeta$, we will consider the following quantity:
$$
h\left(p^{-s^0_1-t(\bs{\beta},\rho)\theta_1},\dots,p^{-s^0_n-t(\bs{\beta},\rho)\theta_n},p^{-c}\right);
$$
and we will prove that for almost all the $t(\bs{\beta},\rho)\in \Xi_{u,\eta}$ (all except a finite number) and for all prime number  $p$ large enough ($p>p_0$ where $p_0$ is an absolute constant) we have:
$$
h\left(p^{-s^0_1-t(\bs{\beta},\rho)\theta_1},\dots,p^{-s^0_n-t(\bs{\beta},\rho)\theta_n},p^{-c}\right)\neq 0.
$$

Write:

\begin{displaymath}
% use packages: array
\begin{array}{lll}
h\left(p^{-s^0_1-t(\bs{\beta},\rho)\theta_1},...,p^{-s^0_n-t(\bs{\beta},\rho)\theta_n},p^{-c}\right) & = & 1+\summ_{k=1}^{r}a_k p^{-\mathbf{s}^0\cdot\bs{\alpha}_{\cdot k}-\bs{\theta}_{(n)}\cdot\bs{\alpha}_{(n) k}\left( \frac{\rho - \sum_{j=1}^{r}\beta_j (\mathbf{s}^0 \cdot\bs{\alpha}_{\cdot j})}{\sum_{j=1}^{r}\beta_j(\bs{\theta}_{(n)}\cdot\bs{\alpha}_{(n) j})}\right)}  \\ 
 & = & 1+\summ_{k=1}^{r}a_k p^{\lambda_{k, \bs{\theta}, \bs{\beta}, \rho}(\bs{\sigma}^0)}
\end{array}
\end{displaymath}
where
$$
\lambda_{k, \bs{\theta}, \bs{\beta}, \rho}(\bs{\sigma}^0) =  -u_{k, \bs{\theta}, \bs{\beta}}(\bs{\sigma}^0) - v_{k, \bs{\theta},\bs{\beta},\rho};
$$
with
$$
u_{k, \bs{\theta},\bs{\beta}}(\bs{\sigma}^0) = \bs{\sigma}^0 \cdot\bs{\alpha}_{\cdot k} - \bs{\theta}_{(n)}\cdot\bs{\alpha}_{(n) k}\frac{\sum_{j=1}^{r}\beta_j (\bs{\sigma}^0 \cdot\bs{\alpha}_{\cdot j})}{\sum_{j=1}^{r} \beta_j(\bs{\theta}_{(n)}\cdot \bs{\alpha}_{(n) j})}
$$
and
\begin{equation}\label{afexpressionvk}
v_{k, \bs{\theta},\bs{\beta},\rho} = \bs{\theta}_{(n)}\cdot\bs{\alpha}_{(n) k}\left(\frac{\rho-i\sum_{j=1}^{r}\beta_j (\bs{\tau}_{(n)}^0, \bs{\alpha}_{(n) j})}{\sum_{j=1}^{r}\beta_j (\bs{\theta}_{(n)}\cdot\bs{\alpha}_{(n) j})}\right)+i\bs{\tau}_{(n)}^0\cdot\bs{\alpha}_{(n) k};
\end{equation}
is independent of $\bs{\sigma}^0$.
Let us precise the dependence of $u_{k, \bs{\theta, \bs{\beta}}}(\bs{\sigma}^0)$ on $\bs{\sigma}^0$.
Indeed, the $(n+1)$-uple $\bs{\sigma}^0$ has here two constraints: $\sigma^0_{n+1}=c$ and $\bs{\sigma}^0\cdot  \bs{\alpha}_{\cdot e}=0$. As in Lemma \ref{aW-noncyclo}, we can consider this  $(n+1)$-uple as a $(n-1)$-uple $\widetilde{\bs{\sigma}}^0=\left(\widetilde{\sigma}^0_1,\dots,\widetilde{\sigma}^0_{n-1}\right)$ without constraint by putting:

\begin{equation*}
\begin{cases}
\displaystyle \sigma^0_{\ell} = \widetilde{\sigma}^0_{\ell} & (\ell \in \{1,\dots,n-1\}), \\
\displaystyle \sigma^0_n = -\frac{1}{\alpha_{n, e}}\left(\sum_{i=1}^{n-1}\alpha_{i, e}\widetilde{\sigma}^0_i + c \ \alpha_{n+1, e}\right) & .
\end{cases}
\end{equation*}

In this way we obtain:

\begin{scriptsize}\begin{displaymath}
% use packages: array
\begin{array}{lll}
u_{k, \bs{\theta, \bs{\beta}}}(\widetilde{\bs{\sigma}}^0) & = &\displaystyle \sum_{i=1}^{n-1}\widetilde{\sigma}^0_i \left(\alpha_{i, k}-\sum_{j=1}^{r}\beta_j\alpha_{i, j} \frac{\bs{\theta}_{(n)}\cdot\bs{\alpha}_{(n) k}}{\sum_{j=1}^{r}\beta_j(\bs{\theta}_{(n)}\cdot  \bs{\alpha}_{(n) j})} -\frac{\alpha_{i, e}}{\alpha_{n, e}}\left(\alpha_{n, k}-\sum_{j=1}^{r}\beta_j\alpha_{n, j} \frac{\bs{\theta}_{(n)}\cdot\bs{\alpha}_{(n) k}}{\sum_{j=1}^{r}\beta_j(\bs{\theta}_{(n)}\cdot  \bs{\alpha}_{(n) j})}\right)\right) \\ 
 &  &\displaystyle +  c \left(\alpha_{n+1, k}-\sum_{j=1}^{r}\beta_j\alpha_{n+1, j} \frac{\bs{\theta}_{(n)}\cdot\bs{\alpha}_{(n) k}}{\sum_{j=1}^{r}\beta_j(\bs{\theta}_{(n)}\cdot  \bs{\alpha}_{(n) j})} -\frac{\alpha_{n+1, e}}{\alpha_{n, e}}\left(\alpha_{n, k}-\sum_{j=1}^{r}\beta_j\alpha_{n, j} \frac{\bs{\theta}_{(n)}\cdot\bs{\alpha}_{(n) k}}{\sum_{j=1}^{r}\beta_j(\bs{\theta}_{(n)}\cdot  \bs{\alpha}_{(n) j})}\right)\right) \\ 
 & = &\displaystyle {u_{k, \bs{\theta, \bs{\beta}}}(\widetilde{\bs{\sigma}}^0)}_{\textrm{vect}} + {u_{k, \bs{\theta}, \bs{\beta}}}_{\textrm{aff}};
\end{array}
\end{displaymath}                 \end{scriptsize}

where

\begin{footnotesize}\begin{equation}\label{aexprukvect}
{u_{k, \bs{\theta, \bs{\beta}}}(\widetilde{\bs{\sigma}}^0)}_{\textrm{vect}}:= \displaystyle \sum_{i=1}^{n-1}\widetilde{\sigma}^0_i \left(\alpha_{i, k}-\sum_{j=1}^{r}\beta_j\alpha_{i, j} \frac{\bs{\theta}_{(n)}\cdot\bs{\alpha}_{(n) k}}{\sum_{j=1}^{r}\beta_j(\bs{\theta}_{(n)}\cdot  \bs{\alpha}_{(n) j})} -\frac{\alpha_{i, e}}{\alpha_{n, e}}\left(\alpha_{n, k}-\sum_{j=1}^{r}\beta_j\alpha_{n, j} \frac{\bs{\theta}_{(n)}\cdot\bs{\alpha}_{(n) k}}{\sum_{j=1}^{r}\beta_j(\bs{\theta}_{(n)}\cdot  \bs{\alpha}_{(n) j})}\right)\right); 
\end{equation}  \end{footnotesize}
and
\begin{footnotesize}
\begin{equation}\label{afexpressionukaff}
{u_{k, \bs{\theta}, \bs{\beta}}}_{\textrm{aff}}:= c \left(\alpha_{n+1, k}-\sum_{j=1}^{r}\beta_j\alpha_{n+1, j} \frac{\bs{\theta}_{(n)}\cdot\bs{\alpha}_{(n) k}}{\sum_{j=1}^{r}\beta_j(\bs{\theta}_{(n)}\cdot  \bs{\alpha}_{(n) j})} -\frac{\alpha_{n+1, e}}{\alpha_{n, e}}\left(\alpha_{n, k}-\sum_{j=1}^{r}\beta_j\alpha_{n, j} \frac{\bs{\theta}_{(n)}\cdot\bs{\alpha}_{(n) k}}{\sum_{j=1}^{r}\beta_j(\bs{\theta}_{(n)}\cdot  \bs{\alpha}_{(n) j})}\right)\right).
\end{equation}  \end{footnotesize}

Then we define the following equivalence relation $\mathcal{R}_{\bs{\beta},\bs{\theta}}$ on the $\alpha_k$

\begin{displaymath}
 % use packages: array
 \begin{array}{lll}
 \displaystyle \bs{\alpha}_{\cdot k} \ \mathcal{R}_{\bs{\beta},\bs{\theta}} \ \bs{\alpha}_{\cdot k'} & \Longleftrightarrow &\displaystyle \textrm{for all} \  \widetilde{\bs{\sigma}}^0 \ \textrm{such that} \ \mathbf{s}_{(n)}^0\in \mathcal{B} \\ 
  &  & {u_{k, \bs{\theta}, \bs{\beta}}}_{\textrm{vect}}(\widetilde{\bs{\sigma}}^0) = {u_{k', \bs{\theta}, \bs{\beta}}}_{\textrm{vect}}(\widetilde{\bs{\sigma}}^0) 
 \end{array}
 \end{displaymath}

Thus $\bs{\alpha}_{\cdot k} \ \mathcal{R}_{\bs{\beta},\bs{\theta}} \ \bs{\alpha}_{\cdot k'}$ if and only if for all $i\in \{1,\dots,n-1\}$:

\begin{small}\begin{equation}\label{afrelationequiv}
% use packages: array
\displaystyle  \left(\alpha_{i, k}-\alpha_{i, k'}\right) - \frac{\alpha_{i, e}}{\alpha_{n, e}}\left(\alpha_{n, k}-\alpha_{n, k'}\right) - \frac{\bs{\theta}_{(n)}\cdot(\bs{\alpha}_{(n) k}-\bs{\alpha}_{(n) k'})}{\sum_{j=1}^{r}\beta_j(\bs{\theta}_{(n)}\cdot  \bs{\alpha}_{(n) j})}\left(\sum_{j=1}^{r}\beta_j\left(\alpha_{i, j}-\frac{\alpha_{i, e}}{\alpha_{n, e}}\alpha_{n, j}\right)\right)=0.
\end{equation}              \end{small}

Notice that
$
\sum_{j=1}^{r}\beta_j\left(\alpha_{i, j}-\frac{\alpha_{i, e}}{\alpha_{n, e}}\alpha_{n, j}\right) = \sum_{j\notin \Lambda_e}\beta_j\left(\alpha_{i, j}-\frac{\alpha_{i, e}}{\alpha_{n, e}}\alpha_{n, j}\right).
$
In addition, it is important to observe that, although the set of $\bs{\beta}\in \Nr$ such that $\gamma(\bs{\beta})\neq 0$ is infinite (since here $h$ is supposed to be non cyclotomic), the set: 
\begin{equation}\label{afensemblefini}
E:=\{\beta_j \mid j\notin \Lambda_e, \gamma(\bs{\beta})\neq 0, \Re(t(\bs{\beta},\rho))\geq 0 \} \ \textrm{is finite.}
\end{equation}

Indeed, since the $t(\bs{\beta},\rho)$ which could cancel the zeros  $t_{m,\bs{\mu},\bs{\sigma}^0}$ are necessarily of positive real part, we have
\begin{equation}\label{afpreellepos}
\Re(t(\bs{\beta},\rho)) = \frac{\Re(\rho)-\sum_{j\notin \Lambda_e}\beta_j(\bs{\sigma}^0\cdot  \bs{\alpha}_{\cdot j})}{\sum_{j=1}^{r}\beta_j(\bs{\theta}\cdot  \bs{\alpha}_{\cdot j})} \geq 0;
\end{equation}

and hence
$
\sum_{j\notin \Lambda_e}\beta_j(\bs{\sigma}^0\cdot  \bs{\alpha}_{\cdot j}) \leq \Re(\rho) <1;
$
which implies that (\ref{afensemblefini}) is a finite set since for all $j\notin \Lambda_e,$ $\bs{\sigma}^0\cdot  \bs{\alpha}_{\cdot j}>0$.
Consequently the quantity $\sum_{j\notin \Lambda_e}\beta_j\left(\alpha_{i, j}-\frac{\alpha_{i, e}}{\alpha_{n, e}}\alpha_{n, j}\right)$ can take only finitely many values when $\bs{\beta}$ moves.

\vspace{0.2cm}

Let us give some precisions concerning the relation $\mathcal{R}_{\bs{\beta},\bs{\theta}}$.

If we assume that $\sum_{j=1}^{r}\beta_j\bs{\alpha}_{(n) j}\notin \mathbb{Q}\bs{\alpha}_{(n) e}$, then there exists $i\in \{1,\dots,n-1\}$ such that $\sum_{j\notin \Lambda_e}\beta_j\left(\alpha_{i, j}-\frac{\alpha_{i, e}}{\alpha_{n, e}}\alpha_{n, j}\right)\neq 0$ and the equality (\ref{afrelationequiv}) is possible only if
\begin{equation}\label{afegalitefiniinfini}
\frac{\bs{\theta}_{(n)}\cdot(\bs{\alpha}_{(n) k}-\bs{\alpha}_{(n) k'})}{\sum_{j=1}^{r}\beta_j(\bs{\theta}_{(n)}\cdot  \bs{\alpha}_{(n) j})} = \frac{\left(\alpha_{i, k}-\alpha_{i, k'}\right) - \frac{\alpha_{i, e}}{\alpha_{n, e}}\left(\alpha_{n, k}-\alpha_{n, k'}\right)}{\sum_{j\notin \Lambda_e}\beta_j\left(\alpha_{i, j}-\frac{\alpha_{i, e}}{\alpha_{n, e}}\alpha_{n, j}\right)}.
\end{equation}

Since the set $\left\{\frac{\left(\alpha_{i, k}-\alpha_{i, k'}\right) - \frac{\alpha_{i, e}}{\alpha_{n, e}}\left(\alpha_{n, k}-\alpha_{n, k'}\right)}{\sum_{j\notin \Lambda_e}\beta_j\left(\alpha_{i, j}-\frac{\alpha_{i, e}}{\alpha_{n, e}}\alpha_{n, j}\right)}\mid \bs{\beta}\in \Nr\right\}$ is finite by (\ref{afensemblefini}) and  $\bs{\theta}_{(n)}\cdot\bs{\alpha}_{(n) j} >0$ for all $j\in \{1,\dots,r\}$, the identity  (\ref{afegalitefiniinfini}) cannot be satisfied for $\lVert\bs{\beta}\lVert$ large enough (i.e. for $\lVert\bs{\beta}\lVert>B_0$ where $B_0$ is an absolute constant) if $\bs{\theta}_{(n)}\cdot(\bs{\alpha}_{(n) k}-\bs{\alpha}_{(n) k'})\neq 0$ or $\left(\alpha_{i, k}-\alpha_{i, k'}\right) - \frac{\alpha_{i, e}}{\alpha_{n, e}}\left(\alpha_{n, k}-\alpha_{n, k'}\right)\neq 0$  because the member on the left is not zero and tends to  $0$ when $\lVert\bs{\beta}\lVert$ tends to infinity.

For $\lVert\bs{\beta}\lVert>B_0$, we have necessarily $\bs{\theta}_{(n)}\cdot(\bs{\alpha}_{(n) k}-\bs{\alpha}_{(n) k'}) = 0$ and $\left(\alpha_{i, k}-\alpha_{i, k'}\right) - \frac{\alpha_{i, e}}{\alpha_{n, e}}\left(\alpha_{n, k}-\alpha_{n, k'}\right) = 0$ for the indexes $i\in \{1,\dots,n-1\}$ such that $\sum_{j\notin \Lambda_e}\beta_j\left(\alpha_{i, j}-\frac{\alpha_{i, e}}{\alpha_{n, e}}\alpha_{n, j}\right)\neq 0$.

For the other indexes $i$ such that $\sum_{j\notin \Lambda_e}\beta_j\left(\alpha_{i, j}-\frac{\alpha_{i, e}}{\alpha_{n, e}}\alpha_{n, j}\right) = 0$, the identity (\ref{afrelationequiv}) also provides $\left(\alpha_{i, k}-\alpha_{i, k'}\right) - \frac{\alpha_{i, e}}{\alpha_{n, e}}\left(\alpha_{n, k}-\alpha_{n, k'}\right) = 0$; thus we obtain   $\bs{\alpha}_{(n) k}-\bs{\alpha}_{(n) k'}\in \mathbb{Q}\bs{\alpha}_{(n) e}$.
And if we write $\bs{\alpha}_{(n) k}-\bs{\alpha}_{(n) k'} = q \bs{\alpha}_{(n) e}$ ($q\in\mathbb{Q}$), the identity $\bs{\theta}_{(n)}\cdot(\bs{\alpha}_{(n) k}-\bs{\alpha}_{(n) k'}) = 0$ gives immediately $q=0$ since $\bs{\theta}_{(n)}\cdot  \bs{\alpha}_{(n) e}\neq 0$; and hence $\bs{\alpha}_{(n) k} = \bs{\alpha}_{(n) k'}$.

Now if $\sum_{j=1}^{r}\beta_j\bs{\alpha}_{(n) j}\in \mathbb{Q}\bs{\alpha}_{(n) e}$, the equality (\ref{afrelationequiv}) becomes:
\begin{equation}\label{afparallalphae}
\forall i \in \{1,\dots,n-1\}, \ \left(\alpha_{i, k}-\alpha_{i, k'}\right)-\frac{\alpha_{i, e}}{\alpha_{n, e}}\left(\alpha_{n, k}-\alpha_{n, k'}\right)=0.
\end{equation}

But then (\ref{afparallalphae}) gives $\bs{\alpha}_{(n) k}-\bs{\alpha}_{(n) k'}\in \mathbb{Q}\bs{\alpha}_{(n) e}$.
Thus, for $\bs{\beta}$ large enough ($\lVert\bs{\beta}\lVert>B_0$), we have:

\begin{displaymath}
% use packages: array
\begin{array}{lll}
\bs{\alpha}_{\cdot k} \ \mathcal{R}_{\bs{\beta},\bs{\theta}} \ \bs{\alpha}_{\cdot k'} & \Longrightarrow  & \displaystyle \begin{cases} 
                                                                             \displaystyle \bs{\alpha}_{(n) k}-\bs{\alpha}_{(n) k'}\in \mathbb{Q}\sum_{j=1}^{r}\beta_j\bs{\alpha}_{(n) j}   \\
										\displaystyle \bs{\alpha}_{(n) k}-\bs{\alpha}_{(n) k'}\in \mathbb{Q}\bs{\alpha}_{(n) e}.
                                                                              \end{cases}
 
\end{array}
\end{displaymath}

We write $[k_0]$ the equivalence class of $k_0$ for the relation $\mathcal{R}_{\bs{\beta},\bs{\theta}}$  and we consider a set $\mathcal{V}$ whose elements are a representative of each equivalence class.

Now if we consider $\sigma^0\longmapsto h\left(p^{-s^0_1-t(\bs{\beta},\rho)\theta_1},...,p^{-s^0_n-t(\bs{\beta},\rho)\theta_n},p^{-c}\right)$ as a function of  $(n-1)$ variables  $f_{\bs{\mu},\bs{\beta},\rho}(\widetilde{\sigma}^0_1,\dots,\widetilde{\sigma}^0_{n-1})$, we can write:

$$
f_{\bs{\mu},\bs{\beta},\rho}(\widetilde{\sigma}^0_1,\dots,\widetilde{\sigma}^0_{n-1}) = 1+\sum_{k_0\in \mathcal{V}}\left(\sum_{k\in [k_0]}a_k p^{-v_{k, \bs{\theta},\bs{\beta},\rho}-{u_{k, \bs{\theta}, \bs{\beta}}}_{\textrm{aff}}}\right)p^{-{u_{k_0, \bs{\theta}, \bs{\beta}}(\widetilde{\bs{\sigma}}^0)}_{\textrm{vect}}};
$$
where the linear forms ${u_{k_0, \bs{\theta}, \bs{\beta}}(\widetilde{\bs{\sigma}}^0)}_{\textrm{vect}}$ are two at a time distinct.

\begin{lemme}\label{abetaefini}
We have for $|\textbf{X}^{\bs{\alpha}_{\cdot j}}| < C$ ($j\in \Lambda_e$) ($C$ being the constant defined in Proposition \ref{expansion}), the following equality:
$$
[h]_{e}(\textbf{X}) = \prod_{\bs{\beta} \in B_{e}}\left(1- \textbf{X}^{\sum_{j \in \Lambda_{e}}\beta_j \bs{\alpha}_{\cdot j}} \right)^{\gamma(\bs{\beta})},
$$
where the right side converges absolutely, and each  $\gamma(\bs{\beta})$ is the integral exponent for the factor indexed by  $\bs{\beta}$ inside the cyclotomic expansion of $h(\textbf{X})$ which is given by Proposition \ref{expansion} page \pageref{expansion}. 
\end{lemme}

\begin{dem}
Firstly put  $d_e = \# \Lambda_e,$ and note the corresponding set   $\Lambda_e = \{j_1 < j_2 < \dots < j_{d_e}\}.$ Then we apply Proposition  \ref{expansion} to the polynomial $[h]_e.$ For the same constant $C$ defined in (\ref{afdist-h-zero}) we have the absolute convergence of the infinite cyclotomic expansion corresponding to $[h]_e (\textbf{X})$ whenever each  $|\textbf{X}^{\bs{\alpha}_{\cdot j}}| < C.$ Notice that here the product of this expansion is taken on all the   $\tilde{\bs{\beta}} \in \mathbb{N}^{d_e} - \{ \textbf{0}\}.$ To each of these $\tilde{\bs{\beta}}$ we can associate a unique  $\bs{\beta} = (\beta_1,\dots, \beta_r) \in B_e$ such that 
 $j_i  \in \Lambda_e$ implies  $\beta_{j_i} = \tilde \beta_i,$ for each $i.$  Consequently,    
$\sum_{j=1}^r  \beta_j \bs{\alpha}_{\cdot j} = \sum_{i=1}^{d_e} \tilde \beta_i \bs{\alpha}_{\cdot j_i}$ if $\bs{\beta} \in B_e.$ Concerning the exponents, we conclude that 
$$\gamma(\tilde{\bs{\beta}}) = \gamma (\bs{\beta}) \qquad \text{for each } \bs{\beta} \in B_e,$$
since the expression of $\gamma (\tilde{\bs{\beta}})$ of Proposition \ref{expansion} page \pageref{expansion} coincides with that of  $\gamma (\bs{\beta})$ 
because the $\bs{\beta} \in B_e$ correspond exactly to the $\tilde{\bs{\beta}}$ with the reindexation that we have just defined. This completes the proof. 
\CQFD

\end{dem}

\begin{rqs}\label{afinitudebeta}
Since we suppose here that $[W_{\bs{\mu},\bs{\sigma}^0}]_e$ and hence that $[h]_e$ is cyclotomic, Lemma \ref{abetaefini} permits to conclude that there is only a finite number of exponents $\gamma(\bs{\beta})\neq 0$ such that $\bs{\beta}\in B_e$.
\end{rqs}

Now recall a classical result that we will use later and which is a consequence of the Weierstrass Preparation Theorem whose proof  can be found in \cite{abh}:

\begin{lemme}\label{athmWeierstrass}
 Let
$f : U \to \mathbb{C}$ be a nonzero holomorphic function  defined on an open set $U\subseteq \mathbb{C}^{n}$. Then the zero locus $f^{-1}(0)$ has   empty interior inside  $\mathbb{C}^n$.
\end{lemme}

Now we can state the last lemma which will allow us to finish the proof of Theorem \ref{afprincipal}:

\begin{lemme}
 Moving $\mathbf{s}_{(n)}^0\in \mathcal{B}\cap \partial \mathbf{W}_c(0)$ if necessary, for almost all  $t(\bs{\beta},\rho)\in \Xi_{u,\eta}$ (i.e. all except a finite number):
$$
h\left(p^{-s^0_1-t(\bs{\beta},\rho)\theta_1},...,p^{-s^0_n-t(\bs{\beta},\rho)\theta_n},p^{-c}\right) \neq 0 \ \textrm{for all prime number} \ p.
$$
\end{lemme}

\begin{dem}

Firstly, since the $\bs{\beta}\in B_e$ such that $\gamma(\bs{\beta})\neq 0$ are in finite number according to Remark  \ref{afinitudebeta}, and since the  $\rho$ such that: 
\begin{displaymath}
\begin{cases}
 t(\bs{\beta},\rho)\in \Xi_{u,\eta} \\
\bs{\beta} \in B_e,
\end{cases}
\end{displaymath}
are necessarily also in finite number, the $t(\bs{\beta},\rho)$ such that $\bs{\beta}\in B_e$ and $\gamma(\bs{\beta})\neq 0$ are in finite number.
So it suffices to consider from now on the $\bs{\beta}\notin B_e$.
Now we want to show that moving $\mathbf{s}_{(n)}^0$ if necessary, the function $f_{\bs{\mu},\bs{\beta},\rho}(\widetilde{\bs{\sigma}}^0)$ is nonzero for all $p$.
Write:

$$
f_{\bs{\mu},\bs{\beta},\rho}(\widetilde{\sigma}^0_1,\dots,\widetilde{\sigma}^0_{n-1}) = 1+\sum_{k_0\in \mathcal{V}}\left(\sum_{k\in [k_0]}a_k  p^{-v_{k, \bs{\theta},\bs{\beta},\rho}-{u_{k, \bs{\theta}, \bs{\beta}}}_{\textrm{aff}}}\right)p^{-{u_{k_0, \bs{\theta}, \bs{\beta}}(\widetilde{\bs{\sigma}}^0)}_{\textrm{vect}}};
$$
where the linear forms  ${u_{k_0, \bs{\theta}, \bs{\beta}}(\widetilde{\bs{\sigma}}^0)}_{\textrm{vect}}$ are two at a time distinct.
Now prove that no ${u_{k_0, \bs{\theta}, \bs{\beta}}(\widetilde{\bs{\sigma}}^0)}_{\textrm{vect}}$ for $k_0\in \mathcal{V}$ is zero for $\beta\notin B_e$ large enough.

So let $k_0\in \mathcal{V}$.
According to the expression of ${u_{k_0, \bs{\theta}, \bs{\beta}}(\widetilde{\bs{\sigma}}^0)}_{\textrm{vect}}$ given in (\ref{aexprukvect}), we have 
\begin{equation}\label{avectequiv}
{u_{k_0, \bs{\theta}, \bs{\beta}}(\widetilde{\bs{\sigma}}^0)}_{\textrm{vect}}=0 \Longleftrightarrow \bs{\alpha}_{(n) k_0}-\sum_{j=1}^{r}\beta_j\bs{\alpha}_{(n) j}\frac{\bs{\theta}_{(n)}\cdot\bs{\alpha}_{(n) k_0}}{\sum_{j=1}^{r}\beta_j(\bs{\theta}_{(n)}\cdot\bs{\alpha}_{(n) j})}\in \mathbb{Q}\bs{\alpha}_{(n) e}.
\end{equation}

Assume that there exists a sequence $(\bs{\beta}_m)_{m\in \mathbb{N}}$ ($\bs{\beta}_{m}\notin B_e$) with $\|\bs{\beta}_m\|\rightarrow +\infty$ when $m\rightarrow +\infty$ so that for all  $m$ there exists $\rho_m$ such that $t(\bs{\beta}_m,\rho_m)\in \Xi_{u,\eta}$ and verifying for all $m$:
\begin{equation}\label{aavantlim}
\bs{\alpha}_{(n) k_0}-\sum_{j=1}^{r}{\beta_{m}}_j\bs{\alpha}_{(n) j}\frac{\bs{\theta}_{(n)}\cdot\bs{\alpha}_{(n) k_0}}{\sum_{j=1}^{r}{\beta_{m}}_j(\bs{\theta}_{(n)}\cdot\bs{\alpha}_{(n) j})}\in \mathbb{Q}\bs{\alpha}_{(n) e}.
\end{equation}

Since (\ref{afensemblefini}) is a finite set we have:
$$
\lim_{m\to +\infty}\frac{\sum_{j=1}^{r}{\beta_{m}}_j\bs{\alpha}_{(n) j}}{\sum_{j=1}^{r}{\beta_{m}}_j(\bs{\theta}_{(n)}\cdot\bs{\alpha}_{(n) j})} = \lim_{m\to +\infty} \frac{\sum_{j\in \Lambda_e}{\beta_{m}}_j\bs{\alpha}_{(n) j}}{\sum_{j\in \Lambda_e}{\beta_{m}}_j(\bs{\theta}_{(n)}\cdot\bs{\alpha}_{(n) j})} = \frac{\bs{\alpha}_{(n) e}}{\bs{\theta}_{(n)}\cdot  \bs{\alpha}_{(n) e}}.
$$

By passing to the limit inside (\ref{aavantlim}) we obtain necessarily $\bs{\alpha}_{(n) k_0}\in \mathbb{Q}\bs{\alpha}_{(n) e}$. 

Consequently, again by  (\ref{aavantlim}), we have that for all $m$:
$$
\sum_{j=1}^{r}{\beta_{m}}_j\bs{\alpha}_{(n) j}\in \mathbb{Q}\bs{\alpha}_{(n) e}.
$$

Hence there exists $q_m\in \mathbb{N}^{*}$ such that $\sum_{j=1}^{r}{\beta_{m}}_j\bs{\alpha}_{(n) j}=q_m \widehat{\bs{\alpha}}_{(n) e}$.

As $c\in \mathbb{Z}\setminus \{0\}$ we have on one hand:

\begin{displaymath}
% use packages: array
\begin{array}{lll}
\bs{\sigma}^0\cdot\left(\sum_{j=1}^{r}{\beta_{m}}_j \bs{\alpha}_{\cdot j}\right) & = &\ds \bs{\sigma}^0_{(n)}\cdot\left(\sum_{j=1}^{r}{\beta_{m}}_j\bs{\alpha}_{(n) j}\right) + c\sum_{j=1}^{r}{\beta_{m}}_j\alpha_{n+1, j}\\ 
 & = &\ds -c q_m \sum_{j=1}^{r}{\beta_{m}}_j\widehat{\alpha}_{n+1, e} + c \sum_{j=1}^{r}{\beta_{m}}_j\alpha_{n+1, j}\in \mathbb{Z}.
\end{array}
\end{displaymath}
On the other hand we know that for all $j$, $\bs{\sigma}^0\cdot  \bs{\alpha}_{\cdot j} \geq 0$ with a strict inequality for  
$j\notin \Lambda_e$. We also know that $\bs{\beta}_{m}\notin B_e$ implies that there exists $j \notin \Lambda_e$ such that 
${\beta_{m}}_j >0$. We deduce that 
$$\bs{\sigma}^0\cdot\left(\sum_{j=1}^{r}{\beta_{m}}_j \bs{\alpha}_{\cdot j}\right) =
\sum_{j=1}^{r}{\beta_{m}}_j \left(\bs{\sigma}^0\cdot  \bs{\alpha}_{\cdot j}\right) = \sum_{j\notin \Lambda_e}{\beta_{m}}_j \left(\bs{\sigma}^0\cdot  \bs{\alpha}_{\cdot j}\right) >0.$$

Consequently we have $\bs{\sigma}^0\cdot\left(\sum_{j=1}^{r}{\beta_{m}}_j \bs{\alpha}_{\cdot j}\right)\geq 1$.
But since $t(\bs{\beta}_m,\rho_m)\in \Xi_{u,\eta}$ we must have: 
$$
0<\Re(t(\beta_m,\rho_m)) = \frac{\Re(\rho_m)-\bs{\sigma}^0\cdot\left(\sum_{j=1}^{r}{\beta_{m}}_j \bs{\alpha}_{\cdot j}\right)}{\sum_{j=1}^{r}{\beta_{m}}_j(\bs{\theta}_{(n)}\cdot\bs{\alpha}_{\cdot j})};
$$
and hence $\Re(\rho_m)>1$; which is impossible and proves that ${u_{k_0, \bs{\theta}, \bs{\beta}}(\widetilde{\bs{\sigma}}^0)}_{\textrm{vect}}$ for $k_0\in \mathcal{V}$ is nonzero for $\beta\notin B_e$ large enough.

Let us consider now $f_{\bs{\mu},\bs{\beta},\rho}(\widetilde{\sigma}^0_1,\dots,\widetilde{\sigma}^0_{n-1})$ and let us prove that $f_{\bs{\mu},\bs{\beta},\rho}(\widetilde{\sigma}^0_1,\dots,\widetilde{\sigma}^0_{n-1})$ is nonzero by moving $\widetilde{\bs{\sigma}}^0$ if necessary.

Firstly, if all the $\sum_{k\in [k_0]}a_k p^{-v_{k, \bs{\theta},\bs{\beta},\rho}-{u_{k, \bs{\theta}, \bs{\beta}}}_{\textrm{aff}}}$ are zero for $k_0\in \mathcal{V}$, then we obtain that 
$f_{\bs{\mu},\bs{\beta},\rho}(\widetilde{\sigma}^0_1,\dots,\widetilde{\sigma}^0_{n-1})$ is a constant function equal to  $1\neq 0$ and satisfies the assertion of the lemma.

Otherwise, there exists at least one $k_0\in \mathcal{V}$ such that $\sum_{k\in [k_0]}a_k p^{-v_{k, \bs{\theta},\bs{\beta},\rho}-{u_{k, \bs{\theta}, \bs{\beta}}}_{\textrm{aff}}}\neq 0$.

Now prove that for all  $\rho$ and $\bs{\beta}\in \Nr$ fixed, the function $\widetilde{\bs{\sigma}}^0\longmapsto f_{\bs{\mu},\bs{\beta},\rho}(\widetilde{\sigma}^0_1,\dots,\widetilde{\sigma}^0_{n-1})$ is nonzero; and this in a way to ensure the fact that its zeros define a thin set of $\mathbb{R}^n$ (i.e. of empty interior).
It suffices for this to consider $\bs{\phi}\in \mathbb{R}^{n-1}$, for example of components $\mathbb{Q}$-linearly independent, so that the  ${u_{k_0, \bs{\theta}, \bs{\beta}}(\bs{\phi})}_{\textrm{vect}}$ are two at a time distinct for $k_0\in \mathcal{V}$.
And we put $
\widetilde{\bs{\sigma}}^0 = t  \bs{\phi}.
$
Since ${u_{k_0, \bs{\theta}, \bs{\beta}}(t\bs{\phi})}_{\textrm{vect}} = t {u_{k_0, \bs{\theta}, \bs{\beta}}(\bs{\phi})}_{\textrm{vect}}$ we obtain:

$$
f_{\bs{\mu},\bs{\beta},\rho}(t\bs{\phi}) = 1+\sum_{k_0\in \mathcal{V}}\left(\sum_{k\in [k_0]}a_k p^{-v_{k, \bs{\theta},\bs{\beta},\rho}-{u_{k, \bs{\theta}, \bs{\beta}}}_{\textrm{aff}}}\right)\exp\left(-t\log(p){u_{k_0, \bs{\theta}, \bs{\beta}}(\bs{\phi})}_{\textrm{vect}}\right).
$$

Then it suffices to use the fact that the functions $\{t\longmapsto \exp\left(-t\log(p){u_{k_0, \bs{\theta}, \bs{\beta}}(\bs{\phi})}_{\textrm{vect}}\right)\}_{k_0\in \mathcal{V}}$ are linearly independent  since the ${u_{k_0, \bs{\theta}, \bs{\beta}}(\bs{\phi})}_{\textrm{vect}}\in \mathbb{R}$ are two at a time distinct; and consequently the function $t\longmapsto f_{\bs{\mu},\bs{\beta},\rho}(t\bs{\phi})$ is nonzero and the function  $\widetilde{\bs{\sigma}}^0\longmapsto f_{\bs{\mu},\bs{\beta},\rho}(\widetilde{\sigma}^0_1,\dots,\widetilde{\sigma}^0_{n-1})$ is also nonzero.
Now it suffices to use Lemma \ref{athmWeierstrass} to deduce that, since $f_{\bs{\mu},\bs{\beta},\rho}(\widetilde{\sigma}^0_1,\dots,\widetilde{\sigma}^0_{n-1})$ is nonzero, the set $f_{\bs{\mu},\bs{\beta},\rho}^{-1}(0)$ is of empty interior inside $\mathbb{C}^{n-1}$ and even inside $\mathbb{R}^{n-1}$ (because any function holomorphic on an open set $U\subseteq \mathbb{C}^n$ and  zero on $U\cap \mathbb{R}^n$ is necessarily zero on  $U$).

Then we put:

$$
\mathfrak{M} := \bigcup_{\bs{\beta}\in \Nr, p,\rho\mid\zeta(\rho)=0}f_{\bs{\mu},\bs{\beta},\rho}^{-1}(0).
$$

This set $\mathfrak{M}$, being a countable union of closed sets of empty interior inside  $\mathbb{R}^{n-1}$, is also of empty interior inside  $\mathbb{R}^{n-1}$ according to Baire's theorem.

As a conclusion, it is possible to choose $\widetilde{\bs{\sigma}}^0\notin \mathfrak{M}$ so that the function $t\longmapsto Z\left(\mathbf{s}_{(n)}^0+t\bs{\theta}_{(n)}\right)$ admits an accumulation of zeros $t_{m,\bs{\mu},\bs{\sigma}^0}$ inside  $\Xi_{u,\eta}$ whithout being cancelled by poles  $t(\bs{\beta},\rho)$; which completes the proof of this lemma and the proof of Theorem  \ref{afprincipal}.
\CQFD

\end{dem}

\section{Proof of Theorem \ref{afprincipal2}.}

To prove Theorem \ref{afprincipal2}, we must localize in the previous proof of Theorem \ref{afprincipal} the use of the hypothesis (H).

In section \ref{afproofthm1} we have considered a vector  $\bs{\alpha}_{\cdot e}$ ($e\in \{1,\dots,r\}$) such that  $\bs{\alpha}_{(n) e}$ determines the polar vector of a face  $\mathcal{F}(\bs{\alpha}_{(n) e})\subseteq\partial \mathbf{W}_c(0)$.
Notice that each face of  $\partial \mathbf{W}_c(0)$ is determined by a polar vector of the form  $\bs{\alpha}_{(n) j}$ for a certain $j$; and if in particular $\bs{\alpha}_{(n) j}\in \mathbb{Q}\bs{\alpha}_{(n) e}$, then the vector $\bs{\alpha}_{\cdot j}$ determines a face of $\partial \mathbf{W}_c(0)$ only if $\bs{\alpha}_{\cdot j}=\bs{\alpha}_{\cdot e}$.
And, having fixed this vector $\bs{\alpha}_{\cdot e}$, we have considered a point  $\mathbf{s}_{(n)}^0$ lying on this face  $\mathcal{F}(\bs{\alpha}_{(n) e})$ (i.e. such that $\bs{\sigma}^0\cdot  \bs{\alpha}_{\cdot e}=0$ and $\bs{\sigma}^0\cdot  \bs{\alpha}_{\cdot j}\geq 0$ for all $j\in \{1,\dots,r\}$).
Then we needed the fact that the vector $\bs{\alpha}_{\cdot e'}$, chosen so that $\bs{\sigma}^0\cdot \bs{\alpha}_{\cdot e'}$ is minimal among the $\bs{\sigma}^0\cdot \bs{\alpha}_{\cdot j_0}>0$ verifying  $\sum_{\{j:\bs{\alpha}_{\cdot j}-\bs{\alpha}_{\cdot j_0}\in \mathbb{Q}\bs{\alpha}_{\cdot e}\}}a_j p^{-i\bs{\tau}^0_{(n)}\cdot \bs{\alpha}_{(n) j}}{c_{\bs{\mu},0}}^{\bs{\theta}_{(n)}\cdot  \bs{\alpha}_{(n) j}}\neq 0$ (see (\ref{afchoixe}) page \pageref{afchoixe}), verifies the following condition:
$$
\bs{\alpha}_{(n) e'}\notin \mathbb{Q} \bs{\alpha}_{(n) e}.
$$
It is only to ensure this property that we have used the hypothesis (H) in the proof of Theorem \ref{afprincipal}.
Obviously, this condition is not a priori satisfied if we do not assume the hypothesis (H).
But the result we want to prove here is weaker than Theorem \ref{afprincipal}.
Indeed, we want to prove the fact that we cannot translate globally the boundary $\partial \mathbf{W}_c(0)$ until $\partial \mathbf{W}_{c}(\delta)$ for all $\delta<0$ if $h$ is not cyclotomic.
And consequently the point $\mathbf{s}_{(n)}^0$ is not constrained to stay in a neighbourhood of a point of $\mathcal{F}(\bs{\alpha}_{(n) e})$ contrary to the previous section  \ref{afproofthm1}.

In addition, we know according to the hypothesis (\ref{afhypnaturelle}) that the set
$$
\mathcal{E}_{e}:= \{\bs{\alpha}_{\cdot j}\mid \bs{\alpha}_{(n) j}\notin \bs{\alpha}_{(n) e}\}\neq \emptyset.
$$

Now let us prove that it is possible to have $\bs{\alpha}_{\cdot e'}\in \mathcal{E}_{e}$ by moving the point $\mathbf{s}_{(n)}^0$ if necessary on the face $\mathcal{F}(\bs{\alpha}_{(n) e})$.

To begin with, let us consider the quantity $\bs{\sigma}^0\cdot  \bs{\alpha}_{\cdot j}$ for all the vectors $\bs{\alpha}_{\cdot j}\notin \mathcal{E}_e$ (i.e. such that $\bs{\alpha}_{(n) j}\in \mathbb{Q}\bs{\alpha}_{(n) e}$). For these $\bs{\alpha}_{\cdot j}$, there exists $q_j\in \mathbb{Q}$ such that $\bs{\alpha}_{(n) j} = q_j \bs{\alpha}_{(n) e}$, and consequently:

\begin{displaymath}
% use packages: array
\begin{array}{lllll}
0 & < & \bs{\sigma}^0\cdot  \bs{\alpha}_{\cdot j} & = & \bs{\sigma}^0_{(n)}\cdot\bs{\alpha}_{(n) j} + c \ \alpha_{n+1, j} \\ 
 &  &  & =  & q_j \bs{\sigma}^0_{(n)}\cdot\bs{\alpha}_{(n) e} + c \ \alpha_{n+1, j} \\ 
 &  &  & = & c\left(\alpha_{n+1, j} - q_j \alpha_{n+1, e}\right) \ \textrm{because} \ \bs{\sigma}^0\cdot  \bs{\alpha}_{\cdot e} = \bs{\sigma}^0_{(n)}\cdot\bs{\alpha}_{(n) e} + c \ \alpha_{n+1, e}=0.  
\end{array}
\end{displaymath}

Thus we observe that $\bs{\sigma}^0\cdot  \bs{\alpha}_{\cdot j}$ for $\alpha_j\notin \mathcal{E}_e$ does not depend on $\bs{\sigma}_{(n)}^0\in \mathcal{F}(\bs{\alpha}_{(n) e})$.
Then we put
$$
\epsilon_0:= \textrm{min}_{\bs{\alpha}_{\cdot j}\notin \mathcal{E}_e} \left(\bs{\sigma}^0\cdot  \bs{\alpha}_{\cdot j}\right)  = \textrm{min}_{\bs{\alpha}_{\cdot j}\notin \mathcal{E}_e}\left(c\left(\alpha_{n+1, j} - q_j \alpha_{n+1, e}\right)\right)>0
$$ 
($\epsilon_0$ does not depend on $\bs{\sigma}_{(n)}^0\in \mathcal{F}(\bs{\alpha}_{(n) e})$).
According to the hypothesis (\ref{afhypnaturelle}) we know that $\partial \mathbf{W}_c(0)$ does not admit only one face because the  $\bs{\alpha}_{\cdot j}\neq \bs{\alpha}_{\cdot e}$ such that $\bs{\alpha}_{(n) j}\in \mathbb{Q}\bs{\alpha}_{(n) e}$ do not define any face of  $\partial \mathbf{W}_c(0)$.
Thus there exists necessarily a vector $\bs{\alpha}_{\cdot j_1}\in \mathcal{E}_e$ such that $\mathcal{F}(\bs{\alpha}_{(n) j_1})$ is a face of $\partial \mathbf{W}_c(0)$ of nonempty intersection with  $\mathcal{F}(\bs{\alpha}_{(n) e})$.
In particular $\mathcal{F}(\bs{\alpha}_{(n) e})\cap \mathcal{F}(\bs{\alpha}_{(n) j_1})$ is also a face of $\partial \mathbf{W}_c(0)$ of dimension strictly inferior.
And the hypothesis (\ref{aeglobal}) ensures the fact that  $\sum_{\{j:\bs{\alpha}_{\cdot j}-\bs{\alpha}_{\cdot j_1}\in \mathbb{Q}\bs{\alpha}_{\cdot e}\}}a_j p^{-i\bs{\tau}^0_{(n)}\cdot \bs{\alpha}_{(n) j}}{c_{\bs{\mu},0}}^{\bs{\theta}_{(n)}\cdot  \bs{\alpha}_{(n) j}}\neq 0$.
So for all $\epsilon>0$, we can find a point  $\mathbf{s}_{(n)}^0\in \mathcal{F}(\bs{\alpha}_{(n) e})$ verifying (\ref{afdcond}) such that
$$
0<\bs{\sigma}^0\cdot\bs{\alpha}_{\cdot j_1}<\epsilon.
$$
And this is in particular true if  $\epsilon<\epsilon_0$.

Now, $e'$ being chosen so that $\bs{\sigma}^0\cdot \bs{\alpha}_{\cdot e'}>0$ is minimal among the $\bs{\sigma}^0\cdot \bs{\alpha}_{\cdot j_0}>0$ verifying  $\sum_{\{j:\bs{\alpha}_{\cdot j}-\bs{\alpha}_{\cdot j_0}\in \mathbb{Q}\bs{\alpha}_{\cdot e}\}}a_j p^{-i\bs{\tau}^0_{(n)}\cdot \bs{\alpha}_{(n) j}}{c_{\bs{\mu},0}}^{\bs{\theta}_{(n)}\cdot  \bs{\alpha}_{(n) j}}\neq 0$, we have:
$$
\bs{\sigma}^0\cdot \bs{\alpha}_{\cdot e'} \leq \bs{\sigma}^0\cdot\bs{\alpha}_{\cdot j_1} <\epsilon_0 = \textrm{min}_{\alpha_j\notin \mathcal{E}_e} \left(\bs{\sigma}^0\cdot  \bs{\alpha}_{\cdot j}\right).
$$
Hence we have necessarily $\bs{\alpha}_{(n) e'}\notin \mathbb{Q}\bs{\alpha}_{(n) e}$; which permits to use again the arguments presented in the section  \ref{afproofthm1} to prove Theorem \ref{afprincipal2}.
\CQFD

\begin{rqs}
We have to notice that the hypothesis (\ref{afhypnaturelle}) is absolutely necessary to ensure as the previous argumentation the existence of a direction  $\bs{\theta}$ in which the zeros or poles of $Z(s_1,\dots,s_n)$ accumulate.
Indeed, consider the following example: $
h(X,Y,Z) = 1 + XY + X^2 Y^2 Z = 1+\textbf{X}^{\bs{\alpha}_{\cdot 1}}+\textbf{X}^{\bs{\alpha}_{\cdot 2}} \in \mathbb{Z}[X,Y,Z];
$
where $\bs{\alpha}_{\cdot 1}=(1,1,0)$ and $\bs{\alpha}_{\cdot 2}=(2,2,1)$;
and the corresponding Euler product (by putting here $c=1$):
$
Z(s_1,s_2) = \prod_{p}h(p^{-s_1},p^{-s_2},p^{-1}). 
$
Let us observe in particular $t\longmapsto Z(\mathbf{s}_{(n)}^0+t\bs{\theta}_{(n)})$ with $\mathbf{s}^0 = (1,-1,1)\in \partial \mathbf{W}(0)\cap \{s_3=1\}$ verifying $\mathbf{s}^0\cdot\bs{\alpha}_{\cdot 1}=0$ and $\bs{\theta} = (\theta_1,\theta_2,0)\in \mathbb{Q}^3\cap\{\theta_3=0\}$.
Then we have $
W_{\bs{\mu}, \bs{\sigma}^0}(X,Y) = 1+X^{\mathbf{s}^0\cdot\bs{\alpha}_{\cdot 1}}Y^{\bs{\theta}_{(n)}\cdot\bs{\alpha}_{\cdot 1}} + X^{\mathbf{s}^0\cdot\bs{\alpha}_{\cdot 2}}Y^{\bs{\theta}_{(n)}\cdot\bs{\alpha}_{\cdot 2}} = 1+Y^{\theta_1+\theta_2} + X Y^{2\theta_1+2\theta_2}
$
and $W_{\bs{\mu}, \bs{\sigma}^0}(p^{-1},p^{-t}) = h\left(p^{-s^0_1-t\theta_1},p^{-s^0_2-t\theta_2},p^{-1}\right)$.
By taking up again the previous notations, the Puiseux series $\Omega_{\bs{\mu},\bs{\sigma}^0}(X)$ verifying $W_{\bs{\mu}, \bs{\sigma}^0}(p^{-1},\Omega_{\bs{\mu},\bs{\sigma}^0}(p^{-1}))=0$ correspond to the branches of $1+T+XT^2$ by putting $T:=Y^{\theta_1+\theta_2}$.
But we check without difficulty that there does not exist any branch which verifies $|\Omega_{\bs{\mu},\bs{\sigma}^0}(X)|<1$ for $|X|$ small.
So there is not any accumulation of zeros of $t\longmapsto h\left(p^{-s^0_1-t\theta_1},p^{-s^0_2-t\theta_2},p^{-1}\right)$ (necessarily of the form $t_{m,\bs{\mu},\bs{\sigma}^0} = -\frac{\log\left(\Omega_{\bs{\mu}, \bs{\sigma}^0}(p^{-1})\right)}{\log(p)}+\frac{2 \pi m i}{\log(p)}$; $m\in\mathbb{Z}, \ p$ prime) of positive real part in a neighbourhood of $\Re(t)=0$.

Moreover, there is not any accumulation of zeros or poles coming from the factor  $t\longmapsto \prod_{\bs{\beta}\in \Nr} \zeta_{M_{\delta}} \left((\mathbf{s}^0_{(n)}+t\bs{\theta}_{(n)},c)\cdot\bs{\alpha}\cdot {}^t \! \bs{\beta} \right)^{-\gamma(\bs{\beta})}$ for all $\delta>0$.
Indeed, these zeros or poles are of the form: $
t(\bs{\beta},\rho)=\frac{\rho-\beta_1\mathbf{s}^0\cdot\bs{\alpha}_{\cdot 1}-\beta_2\langle s^0,\bs{\alpha}_{\cdot 2}\rangle}{\beta_1\bs{\theta}_{(n)}\cdot\bs{\alpha}_{\cdot 1}+\beta_2\bs{\theta}_{(n)}\cdot\bs{\alpha}_{\cdot 2}} = \frac{\rho-\beta_2}{(\beta_1+2\beta_2)(\theta_1+\theta_2)},
$
where $\bs{\beta} =(\beta_1,\beta_2)\in \mathbb{N}^2$ and $\rho$ is a zero or a pole of the Riemann zeta function.
So if $\beta_2>0$ we will have $\Re(t(\bs{\beta},\rho)) \leq 0$; hence the $t(\bs{\beta},\rho)$ of positive real part are such that $\beta_2=0$.
But we know (see Lemma \ref{abetaefini}) that these zeros or poles  come only from a finite number of $\zeta$-factors of  $\prod_{\bs{\beta}\in \Nr} \zeta_{M_{\delta}} \left((\mathbf{s}^0_{(n)}+t\bs{\theta}_{(n)},c)\cdot\bs{\alpha}\cdot {}^t \! \bs{\beta} \right)^{-\gamma(\bs{\beta})}$ (which correspond to the cyclotomic factorization of the cyclotomic polynomial  $1+Y^{\theta_1+\theta_2}$); and consequently these zeros or poles are isolated and do not accumulate in a neighbourhood on the right of  $\partial \mathbf{W}(0)$.

\end{rqs}

\section{Proof of Theorem \ref{afapplication}.}

To establish Theorem $\ref{afapplication}$, it suffices to rewrite the Igusa zeta function under the form of a pseudo-uniform Euler product associated to a certain polynomial; and we will check that this polynomial satisfies the conditions of Theorem  $\ref{afprincipal}$.
For $\sigma_i>1$ ($i=1,\dots,n$), let us write:
\begin{displaymath}
% use packages: array
\begin{array}{lll}
Z^{\textrm{ring}}(s_1,\dots,s_n; \mathbb{Z}[T,T^{-1}]) & = & \summ_{m_1,\dots,m_n \geq 1} \frac{\varphi(m_1\cdots m_n)}{m_1^{s_1} \cdots m_n^{s_n}}   \\ 
 & = & \displaystyle\prod_{p}\left(\summ_{\bs{\nu} \in \mathbb{N}^n}\frac{\varphi(p^{\lVert\bs{\nu}\lVert})}{p^{\mathbf{s}\cdot {}^t \! \bs{\nu}}}\right) \\ 
 & = & \displaystyle\prod_{p}\left(1+ \summ_{\bs{\nu} \in \mathbb{N}^n\setminus \{\textbf{0}\}}\frac{\left(p^{\lVert\bs{\nu}\lVert}-p^{\lVert\bs{\nu}\lVert-1}\right)}{p^{\mathbf{s}\cdot {}^t \! \bs{\nu}}}\right) \\
 & = & \displaystyle\prod_{p}\left(1+ \left(\summ_{\bs{\nu} \in \mathbb{N}^n\setminus \{\textbf{0}\}}\frac{p^{\lVert\bs{\nu}\lVert}}{p^{\mathbf{s}\cdot {}^t \! \bs{\nu}}}\right)\left(1-\frac{1}{p}\right)\right).
\end{array}
\end{displaymath}
But for $\sigma_i>2$ ($i=1,\dots,n$), we have:
\begin{displaymath}
% use packages: array
\begin{array}{lll}
\summ_{\bs{\nu} \in \mathbb{N}^n}\frac{p^{\lVert\bs{\nu}\lVert}}{p^{\mathbf{s}\cdot {}^t \! \bs{\nu}}} & = & \summ_{\bs{\nu} \in \mathbb{N}^n}\frac{1}{\prod_{i=1}^{n}p^{\nu_i(s_i-1)}} \\ 
 & = & \summ_{\bs{\nu} \in \mathbb{N}^n}\prod_{i=1}^{n}\left(\frac{1}{p^{(s_i-1)}}\right)^{\nu_i} \\ 
 & = & \displaystyle\prod_{i=1}^{n}\left(\sum_{\nu_i=0}^{+\infty}\left(\frac{1}{p^{(s_i-1)}}\right)^{\nu_i}\right) \\
 & = & \displaystyle\prod_{i=1}^{n}\left(1-\frac{1}{p^{(s_i-1)}}\right)^{-1}.
\end{array}
\end{displaymath}
Consequently, for $\sigma_i>2$ ($i=1,\dots,n$) we obtain:
\begin{displaymath}
% use packages: array
\begin{array}{lll}
Z^{\textrm{ring}}(s_1,\dots,s_n; \mathbb{Z}[T,T^{-1}]) & = & \displaystyle\prod_{p}\left(1+\left(\prod_{i=1}^{n}\left(1-\frac{1}{p^{(s_i-1)}}\right)^{-1}-1\right)\left(1-\frac{1}{p}\right)\right) \\ 
 & = & \displaystyle\prod_{p}\prod_{i=1}^{n}\left(1-\frac{1}{p^{(s_i-1)}}\right)^{-1} \\   &   & \displaystyle\left(\prod_{i=1}^{n}\left(1-\frac{1}{p^{(s_i-1)}}\right) + 
  \left(1-\prod_{i=1}^{n}\left(1-\frac{1}{p^{(s_i-1)}}\right)\right)\left(1-\frac{1}{p}\right)\right) \\ 
 & = & \displaystyle\prod_{i=1}^{n}\zeta\left(s_i-1\right)\prod_{p}\left(1-\frac{1}{p}+\frac{1}{p}\prod_{i=1}^{n}\left(1-\frac{1}{p^{(s_i-1)}}\right)\right) \\
 & = & \displaystyle\prod_{i=1}^{n}\zeta\left(s_i-1\right)\prod_{p}\left(1+\sum_{k=1}^{n}\sum_{\stackrel{I\subseteq \{1,\dots,n\},}{ \#I=k}}\frac{(-1)^{k}}{p^{\left(\sum_{\ell \in I}s_{\ell}\right)-k+1}}\right).
\end{array}
\end{displaymath}
Since the finite product of zeta functions $\displaystyle\prod_{i=1}^{n}\zeta\left(s_i-1\right)$ is meromorphic to whole  $\mathbb{C}^n$, it suffices to determine the maximal domain of meromorphy of the product:
$$
\mathbf{s}\longmapsto\displaystyle\prod_{p}\left(1+\sum_{k=1}^{n}\sum_{\stackrel{I\subseteq \{1,\dots,n\},}{ \#I=k}}\frac{(-1)^{k}}{p^{\left(\sum_{\ell \in I}s_{\ell}\right)-k+1}}\right).
$$
By establishing the change of variable $\textbf{w}=\mathbf{s}-\textbf{1} = (s_1-1,\dots,s_n-1)$, we are led to consider the product:
$$
\textbf{w}\longmapsto\displaystyle\prod_{p}\left(1+\sum_{k=1}^{n}\sum_{\stackrel{I\subseteq \{1,\dots,n\},}{ \#I=k}}\frac{(-1)^{k}}{p^{\left(\sum_{\ell \in I}w_l\right)+1}}\right);
$$
which is equal to:
$$
\displaystyle\prod_{p}h(p^{-w_1},\dots,p^{-w_n},p^{-1}),
$$
with
$$
h(X_1,\dots,X_{n+1}) = 1+\sum_{I\subseteq \{1,\dots,n\}} (-1)^{\#I}\textbf{X}^{\bs{\alpha}_{\cdot I}},
$$
by putting for all  $I\subseteq\{1,\dots,n\},$ $\alpha_{n+1, I}=1$ and for $\ell \in \{1,\dots,n\}$:

\begin{equation*}
 \begin{cases}
 \alpha_{\ell, I} = 1 \ \textrm{if} \ \ell \in I \\
 \alpha_{\ell, I} = 0 \ \textrm{otherwise}.
 \end{cases}
 \end{equation*}

To finish, we can easily check that  $h$ satisfies the conditions of Theorem  $\ref{afprincipal}$; which completes the proof of Theorem $\ref{afapplication}$.
\CQFD

\vspace{0.3cm}

\textbf{ACKNOWLEDGEMENTS:}

This article is the result of a work that the author has done during his thesis. The author wants particularly to thank his advisor Driss Essouabri for his relevant remarks and his careful reading of this paper. He also expresses his gratitude towards Ben Lichtin for his precious advice which have improved the presentation of this paper.

\end{document}